\documentclass[12pt,reqno,psfig,openbib]{amsart}

\usepackage[pagebackref=true, colorlinks=true, citecolor=blue]{hyperref}
\usepackage{amssymb,amsmath,graphicx,amsfonts,euscript}
\usepackage{fancyhdr}
\usepackage{epsfig}
\usepackage{setspace}
\usepackage{dsfont}
\usepackage{graphics}
\usepackage{subcaption}
\usepackage[all]{xy}
\usepackage{amsmath}
\usepackage{enumitem}
\graphicspath{ {./figures/} }

\usepackage{caption}
\usepackage{subcaption}

\usepackage{soul}

\usepackage{float}

 \usepackage{pst-grad} %
 \usepackage{pst-plot} %

\numberwithin{equation}{section}
\numberwithin{figure}{section}

 \usepackage[left=1.25in,top=1.15in,right=1.25in]{geometry}
\usepackage{verbatim}

\newcommand{\spacer}{\vspace{1.5mm}}

\date{\today}

\renewcommand{\epsilon}{\varepsilon}
\newcommand{\ep}{\varepsilon}

\newcommand{\pa}{\partial}
\newcommand{\tht}{\theta}

\newcommand{\pri}{\prime}

\newcommand{\kap}{\kappa}
\newcommand{\blds}{\boldsymbol}


\begin{document}

\title[
Finite Volume Neural Networks for singular  perturbations
]
{Finite volume neural networks for singularly perturbed boundary value problems}

\author[  
G.-M. Gie, 
Junghwa Kim
]
{
Gung-Min Gie$^{1, 2}$, 
Junghwa Kim$^1$
}
\address{$^1$ Department of Mathematics, University of Louisville, Louisville, KY 40292}
\address{$^2$ Research Institute of Mathematics, Seoul National University, Seoul, Korea}
\email{gungmin.gie@louisville.edu}
\email{junghwa.kim@louisville.edu}

\begin{abstract}

We develop a singular-layer enhanced Finite Volume Neural Network (sl‑FVNN) methodology for solving singularly perturbed boundary value problems.
The proposed framework combines flux-conserving finite volume residuals with analytically derived boundary layer correctors, enabling stable and efficient approximation of sharp boundary variations induced by small perturbation parameters. 
Numerical experiments over a broad range of small perturbation parameters ($10^{-4}$ to $10^{-12}$), demonstrate that the sl‑FVNN yields uniformly accurate predicted solutions and significantly improved convergence behavior compared to conventional FVNN methods.    
Further comparisons with sl-PINNs from \cite{GHJ} show that our sl‑FVNNs attain comparable accuracy for reaction–diffusion and Burgers equations, while achieving significantly improved convergence for convection-dominated nonlinear problems. 
These results highlight the effectiveness of incorporating asymptotic layer structures into finite volume based neural networks, offering a practical and scalable framework for multiscale PDE learning in stiff regimes.
\end{abstract}

\maketitle

\tableofcontents

\section{Introduction}\label{S. Intro} 
Singularly perturbed problems arise when a small parameter multiplies the terms with the highest order derivative, e.g., reaction–diffusion or convection–diffusion equations with small diffusivity, and 
the Burgers or Navier-Stokes equations at a small viscosity.  
As the perturbation parameter (small diffusivity or viscosity) tends to zero, the solution often develops sharp boundary, interior, or corner layer regions where the solution changes rapidly within the thin regions.
It is well-known that these features pose significant challenges for both analysis and computation; see, e.g., \cite{Book16, HK82,  OM08, SK87}.  
From the computational point of view, 
{
standard numerical methods on uniform meshes often require highly refined meshes to resolve such layers}, while  
standard Physics Informed Neural Networks (PINNs) suffer from an imbalance between the small parameter and the neural network regularization scale.  
As a result, both approaches often exhibit nonuniform convergence and large local errors near the singular layers.

Physics-Informed Neural Networks (PINNs) have recently emerged as a powerful framework for solving PDEs by embedding the governing equations and boundary conditions into a neural network loss function using automatic differentiation. 
This approach enables the machine learning prediction to approximate smooth solution representations by implicitly enforcing physical laws. 
The PINNs have shown promising results across diverse PDE problems, including fluid dynamics, elasticity, and reaction–diffusion systems; see, e.g., \cite{
KKLPWY21, 
KZK19, 
LMMK21, 
LPYWVJ21, 
LLF98, 
MZK20, 
RPK19, 
GLK19, 
WTP21, 
YMK21}. 
However, 
standard PINNs struggle to approximate any singularly perturbed problem because the presence of thin layers induces stiffness that challenges both training stability and convergence to the solution; see, e.g., \cite{
ACD, 
PINN_bl01, 
Tony9, 
pinn_lim, 
FSF22, 
FT20, 
PINN_bl02}.

To address this issue, several recent approaches have incorporated {\it singular perturbation analysis} into PINN frameworks. One such development is the \emph{singular-layer PINNs (sl-PINNs)}, which enrich the network architecture by adding   
analytic information of the layer structures via auxiliary corrector functions \cite{CGHJ_Burgers, CGHJ_PPF, GHJL, GHJM, GHJ, JKN25}. 
These correctors explicitly represent the dominant layer component, while the neural network learns the smoother residual component, leading to improved convergence and robustness in stiff regimes.

The Finite Volume (FV) method is one of the most widely used numerical schemes for PDEs, 
particularly in fluid mechanics and conservation laws; see, e.g., \cite{EGH06, GT10, GT15, LeV02, MS91} for general theory and \cite{BGHL13, GJN19} for numerical applications. 
The main idea of FV method is to integrate the governing equation over each control volume and apply the divergence theorem, 
and hence enforcing local conservation of fluxes across interfaces.  
This feature provides superior robustness and physical precision, especially for advection-dominated or nonlinear systems where conservation properties play an important role.  
The Finite Volume framework also allows flexible implementation on complex geometries and unstructured meshes. Despite its strengths, FV methods  in general require careful discretization, high mesh resolution near singular layers to capture steep gradients accurately; see., e.g., \cite{GJL1, JT_FV_09, JT10}.

In this article, we combine the strengths of both PINNs and FV methods by introducing a hybrid framework called the \textbf{Finite Volume (Physics Informed) Neural Network (FVNN)}, where the PINN loss is formulated based on flux conservative, Finite Volume residuals rather than differential residuals; see some related earlier works, e.g., \cite{FVPINN1, FVPINN2, FVPINN3}.  
This FVNN reduces the required derivative order, improves numerical stability, and ensures local conservation of fluxes at the discrete level. Building further, we propose a new method, called the \textbf{Singular-Layer Finite Volume Neural Network (sl-FVNN)}, 
which incorporates analytical boundary layer correctors of singular layers into the FVNN architecture. 
These correctors capture the stiff component of the solution (the boundary or interior layer), allowing the neural network to focus on learning the smooth remainder.

The key characteristics and features of the proposed approach are:
\begin{itemize}
    \item[(i)] A Finite Volume Neural Network formulation that preserves local fluxes and reduces the sensitivity to small perturbation parameter by using lower order derivatives in the residual;
    \item[(ii)] Integration of singular layer correctors into the network structure, enhancing accuracy and convergence for singularly perturbed problems;
    \item[(iii)] Demonstrating the effectiveness via numerical experiments and showing uniform accuracy and consistency even in small perturbation parameter regimes.
\end{itemize}

Overall, the sl-FVNN methods provide a stable, flux-conserving, and asymptotically informed learning framework that can resolve sharp boundary features without relying on dense meshing or sacrificing stability.

{
The present work focuses on one-dimensional benchmark problems for which the layer locations and leading-order correctors are explicitly available. Extensions to multidimensional problems require additional asymptotic information on the layer geometry and are planned for future work.
}

\smallskip

The remainder of this paper is organized as follows. 

In Section \ref{S.1D}, we introduce the class of singularly perturbed PDE models including boundary layers that causes computational challenges. For these models, we recall the FVNNs and outline the core architecture of our new sl-FVNNs where we embed the corresponding layer correctors into the FVNN framework.  
Section \ref{S.Comp_1d} demonstrates the effectiveness of the proposed sl-FVNN methods on singularly perturbed 
reaction diffusion equation, 
nonlinear convection diffusion equation, 
Burgers equation. 
We report numerical convergence and robustness with respect to the small perturbation parameter$\ep$ as well as  
comparisons with standard FVNNs, and PINNs and sl-PINNs studied in \cite{GHJ}. Finally, Section \ref{S.Con} contains concluding remarks and future directions.

\section{Finite Volume Neural Networks (FVNNs) for 1D Problems}\label{S.1D}

We consider the 1D singularly perturbed boundary value problem in the form,
\begin{equation}\label{e:1d_gen}
    \left\{
    \begin{array}{rl}
        \spacer
        - \ep d_x^2 u^\ep  + g(u^\ep)   = f(x), & 0 < x <1, \\
        u^\ep(0) = u_0,  & u^\ep(1) = u_1, 
    \end{array}
    \right.
\end{equation}
where $\ep>0$ is the small perturbation parameter, and  
$f = f(x)$ and $u_i \in \mathbb{R}$, $i=0, 1$, are given data. 
The function $g(u^\ep)$ is a linear or nonlinear function in $u^\ep$ that appears in  lower-order derivatives, as in the examples (\ref{e:rd_1D_eqn}), (\ref{e:cd_1d}), and (\ref{e:Burgers}) below.   
\smallskip

To construct the predicted solution to (\ref{e:1d_gen}) via {\it Physics Informed Neural Networks} (PINNs) or {\it Finite Volume (Physics Informed) Neural Networks} (FVNNs), 
we introduce an $L$-layer feed-forward Neural Network (NN) 
from $\mathbb{R}^{N_0}$ to $\mathbb{R}^{N_L}$: 

Define recursively $\mathcal{N}^l$, $0 \leq l \leq L$, where
\begin{equation}\label{e:NN}
\begin{array}{rl}
    \text{input layer: } & \mathcal{N}^{0}(\blds{x}) = \blds{x} \in \mathbb{R}^{N_0}, \\
    \text{hidden layers: } & \mathcal{N}^{l}(\blds{x}) = \sigma(\blds{W}^{l} \mathcal{N}^{l-1}(\blds{x}) + \blds{b}^{l}) \in \mathbb{R}^{N_{l}}, \quad 1\leq l \leq L-1, \\
    \text{output layer: } & \mathcal{N}^{L}(\blds{x}) = \blds{W}^{L} \mathcal{N}^{L-1}(\blds{x}) + \blds{b}^{L} \in \mathbb{R}^{N_{L}},
\end{array}
\end{equation}
where $\blds{x}$ and $\mathcal{N}^{L}(\blds{x})$ are the input and output vectors, and 
$\blds{W}^{l} \in \mathbb{R}^{N_{l} \times N_{l-1}}$ and $\blds{b}^{l} \in \mathbb{R}^{N_{l}}$, $l=1,\ldots,L$, 
are the (unknown) weights and biases respectively. 
Here, $N^{l}$, $l=0,\ldots,L$, is the number of neurons in the $l$-th layer, and $\sigma$ is the activation function. 
We set $\blds{\tht}=\{ (\blds{W}^{l}, \blds{b}^{l}) \}_{l=1}^{L}$ as the collection of all the parameters, weights and biases. 
\smallskip

The \textit{Physics-Informed Neural Network} (PINN) predicted solution to (\ref{e:1d_gen}) is obtained by the NN structure (\ref{e:NN}) with input $x \in \mathbb{R}$ and output $ \mathcal{N}^{L} = \hat{u} (x ;\blds{\tht}) \in \mathbb{R}$. 
This PINN predicted solution $\hat{u}$ is trained by minimizing the loss function, 
\begin{equation}\label{e:loss_PINN}
    \begin{split}
        \mathcal{L}_{PINN}(\blds{\tht}; \mathcal{T}) 
        &= \frac{1}{M} \sum_{x_i \in [0, 1]} \left| 
        - \ep d_x^2 \hat{u}   + g(\hat{u} )   - f  \right|^{2} 
        +
        | \hat{u} (0) - u_0 |^2  +
        | \hat{u} (1) - u_1 |^2,
    \end{split}
\end{equation}
where $\mathcal{T} = \{ x_i \}_{i=1}^{M} \subset [0, 1]$ is a set of training data points. 
By minimizing the loss function, we obtain the PINN prediction $\hat{u}$ to (\ref{e:1d_gen}).  
\smallskip

Concerning the {\it Finite Volume Neural Network} (FVNN) framework for (\ref{e:1d_gen}), 
we use the same set $\mathcal{T} = \{x_i\}_{i=1}^M \subset [0,1]$ with uniformly distributed points and define the control volumes, 
\begin{equation}\label{e:control_volume_1d}
    K_i := (x_{i-1}, x_{i+1}), \quad 2 \leq i \leq M-1,  
\end{equation}
i.e.,    
$    x_i \text{ is the center of } K_i$.

By integrating the equation (\ref{e:1d_gen}) over each $K_i$ 
and applying the divergence theorem, we write 
\begin{equation}\label{e:1d_gen_w}
        - \ep d_x u^\ep(x_{i+1}) 
        + \ep d_x u^\ep(x_{i-1})        
        + \int_{K_i} g(u^\ep) \, dx
    =        
        \int_{K_i} f \, dx; 
\end{equation}
{
hence the local flux of (\ref{e:1d_gen}) is conserved in each $K_i$. }

To define the FVNN prediction $\hat{u}(x; \boldsymbol{\theta})$ of (\ref{e:1d_gen}),  
we minimize the Finite Volume loss function, including the local flux  (\ref{e:1d_gen_w}), defined by  
\begin{equation}\label{e:loss_FVNN}
    \begin{split}
        \mathcal{L}_{FV}(\blds{\tht}; \mathcal{T}) 
        =
        &
        \frac{1}{M-2} \sum_{i=2}^{M-1} 
        \Big| 
        - \ep d_x {\hat{u}} (x_{i+1}) 
        + \ep d_x {\hat{u}} (x_{i-1})
        + g(\hat{u} (x_i))  |K_i| 
        - f(x_i) |K_i|  
        \Big|^{2} \\
        &
        +
        | \hat{u} (0) - u_0 |^2  +
        | \hat{u} (1) - u_1 |^2.
    \end{split}
\end{equation}
The flux associated with the lower-order term $g(\cdot)$ in (\ref{e:1d_gen_w}) and (\ref{e:loss_FVNN}) 
is formulated in a specific way depending on the explicit expression of $g(\cdot)$; see the equations  
(\ref{e:loss_PINNs_mod_rea-dif}), 
(\ref{e:loss_FVNN_cd_1d}), 
(\ref{e:loss_sl_FVNN_cd_1d}), 
(\ref{e:loss_FVNN_Burgers}), and 
(\ref{e:loss_sl_FVNN_Burgers}). 

As we will show in Section \ref{S.Comp_1d} below, 
both PINNs and FVNNs produce relatively accurate solutions to (\ref{e:1d_gen}) when the perturbation parameter $\ep$ is not too small. 
However, neither method approximates the solution well in the small $\ep$ regime; see the tables   \ref{table1_rd},  \ref{table1_cd}, and   \ref{table1_Burgers}, and the other tables and graphs. 

To resolve these issues, we propose a semi-analytic FVNN approach called the singular-layer FVNN
(sl-FVNN), tailored for singular perturbation problems with boundary layers. 
This approach incorporates the boundary layer correctors (analytic approximations of solutions' singular behavior) into the network.
These correctors capture the steep component of the solution, while the neural network learns the smooth remainder.

Combined with the finite volume residual formulation, our novel strategy of sl-FVNNs preserves the local flux and enables
accurate resolution even for extremely small parameter $\ep$. 
Numerical experiments in Section \ref{S.Comp_1d}
confirm that our sl-FVNN framework achieves uniform accuracy with physical consistency across
singular perturbation regimes.

All the examples used in this article are those investigated in the earlier paper \cite{GHJ} where the semi-analytic version of PINNs (sl-PINNs as named afterwards in following papers \cite{CGHJ_Burgers, CGHJ_PPF, GHJL, GHJM, GHJ, JKN25}) 
were first introduced. 
We briefly review those results in Section 2.4 of \cite{GHJ} and compare them to the results of FVNNs and the new sl-FVNNs introduced in this article.

\subsection{Reaction Diffusion Equation (RD eqn)}\label{S.rd_1D}
With $g(u^\ep) = u^\ep$ and $u_0= u_1 = 0$ in (\ref{e:1d_gen}), we consider  
\begin{equation}\label{e:rd_1D_eqn}
    \left\{
    \begin{array}{rl}
        - \ep d_x^2  u^\ep  + u^\ep   = f(x), & 0 < x <1, \\
        u^\ep = 0, & x=0, 1.
    \end{array}
    \right.
\end{equation}

Our main objective  here is 
to construct the semi-analytic sl-FVNN for the singular perturbation problem above, and 
compare its performance with the FVNN as well as 
the PINN and sl-PINN studied in \cite{GHJ}.  
The FVNN loss function in (\ref{e:loss_FVNN}) for this example 
is explicitly given by 
\begin{equation}\label{e:loss_FVNN_rd_1d}
    \begin{split}
        \mathcal{L}_{FV}(\blds{\tht}; \mathcal{T}) 
        =
        &
        \frac{1}{M-2} \sum_{i=2}^{M-1} 
        \Big| 
        - \ep d_x {\hat{u}} (x_{i+1}) 
        + \ep d_x {\hat{u}} (x_{i-1})
        + \hat{u} (x_i)  |K_i| 
        - f(x_i) |K_i|  
        \Big|^{2} \\
        &
        +
        | \hat{u} (0)  |^2  +
        | \hat{u} (1)  |^2.
    \end{split}
\end{equation}
\smallskip

Now, to construct our new singular layer FVNN (sl-FVNN) for the singularly perturbed problem (\ref{e:rd_1D_eqn}), 
we briefly recall the associated boundary layer analysis; see, e.g., \cite{Book16} for details:  

Formally setting $\ep = 0$ in (\ref{e:rd_1D_eqn}), 
we find that
the corresponding limit of (\ref{e:rd_1D_eqn}) at the zero 
diffusivity is
\begin{equation}\label{e:rea_dif_eqn_limit}
    u^0   = f.
\end{equation}

Employing the method of matching asymptotics {in, e.g., \cite{Book16}},  
we {first} notice that the size of boundary layers of this problem (\ref{e:rd_1D_eqn}) near $x = 0$ and $x = 1$ is of order $\ep^{1/2}$. 
{
Considering the difference between the equations (\ref{e:rd_1D_eqn}) and (\ref{e:rea_dif_eqn_limit}), 
we notice that 
\begin{equation*}\label{e:REV_1} 
    \left\{
        \begin{array}{rl}
            \spacer
            -\ep d^2_x (u^\ep - u^0) 
            + (u^\ep - u^0)  =  \ep d^2_x u^0, 
            &
            0 < x < 1, \\
            (u^\ep - u^0) = - u^0 , 
            &
            x = 0, 1.
        \end{array}
    \right.
\end{equation*}
Thus, after dropping the lower order term of order $\ep$ on the right-hand side, 
we write}
the equation of the boundary layer corrector 
$\varphi  \simeq u^\ep - u^0$ as 
\begin{equation}\label{e:cor_eqn_rea-dif}
    \left\{
        \begin{array}{rl}
            \spacer
            -\ep d^2_x \varphi  + \varphi  = 0, 
            &
            0 < x < 1, \\
            \varphi= - u^0 , 
            &
            x = 0, 1.
        \end{array}
    \right.
\end{equation} 

{
Solving (\ref{e:cor_eqn_rea-dif}) directly, 
we obtain 
$$
\varphi(x) = \frac{-u^0(0)+e^{-1/\sqrt{\varepsilon}}u^0(1)} {1-e^{-2/\sqrt{\varepsilon}}} e^{-x/\sqrt{\varepsilon}} + \frac{e^{-1/\sqrt{\varepsilon}}u^0(0)-u^0(1)} {1-e^{-2/\sqrt{\varepsilon}}} e^{-(1-x)/\sqrt{\varepsilon}}.  
$$
Since \(e^{-1/\sqrt{\varepsilon}}\) is exponentially small as \(\varepsilon\to0\), 
we write the explicit expression of the corrector $\varphi$  as 
\begin{equation}\label{e:cor__rea-dif}
    \varphi(x)
        =
            -u^0(0) \, e^{-x/\sqrt{\ep}}
            -u^0(1) \, e^{-(1-x)/\sqrt{\ep}}
            + e.s.t., 
\end{equation}
where \({e.s.t.}\) denotes terms that are exponentially small as \(\varepsilon\to0\).
}

The corrector $\varphi(x)$ (fast decaying part of $u^\ep$ near $x = 0, 1$) is in fact 
the sum of 
two exponentially decaying functions from $x = 0, 1$, 
with respect to the stretched variables $x/\sqrt{\ep}$ and 
$(1-x)/\sqrt{\ep}$, up to an exponentially small term denoted by $e.s.t$.
\smallskip

Following the asymptotic analysis as in \cite{Book16}, 
one can perform the energy estimates on the difference 
$u^\ep - (u^0 + \varphi)$ and 
find that 
\begin{equation}\label{e:con__rea-dif}
        \|
            u^\ep - (u^0 + \varphi)
        \|_{L^2((0, 1))}
            \leq
                \kap \ep,  \qquad
        \|
            u^\ep -  u^0  
        \|_{L^2((0, 1))}
            \leq
                \kap \ep^{\frac{1}{4}}, 
\end{equation}
for a constant $\kap>0$ independent of $\ep$. 

{
We emphasize that the estimates in (\ref{e:con__rea-dif}) are not intended to provide an approximation error bound or a convergence guarantee for the neural network training. 
Rather, they justify the decomposition of $u^\ep$ 
into a singular boundary layer component and a smoother remainder. 
In particular, 
the asymptotic expansion   
$u^0 + \varphi$ 
satisfies the original singularly perturbed equation up to a small remainder as $\ep$ tends to $0$. Therefore, embedding the corrector  into the FVNN ansatz removes the dominant singular layer contribution from the residual and helps prevent the corresponding PINN/FVNN loss from blowing up in the small-$\ep$ regime.
}

We infer from (\ref{e:con__rea-dif}) that 
the corrector $\varphi$ 
exhibits the singular behavior of $u^\ep$ at a small diffusivity $\ep$, 
and 
that 
$u^\ep$ is decomposed into the sum of 
fast (decaying part) $\varphi$, 
and the slow part $u^0$. 
\smallskip

Using the profile of $\varphi$ in (\ref{e:cor__rea-dif}), 
we propose    
our new {\it singular layer FVNN} (sl-FVNN) for  (\ref{e:rd_1D_eqn}) in the form,
\begin{equation}\label{e:app_sol_enriched_rea-def}
    \widetilde{u}(x ; \, {\blds \theta}) 
        =
                \hat{u}({ {x}}; \, {\blds \theta})
            -
                \hat{u}({ {0}}; \, {\blds \theta}) \,
                e^{-x/\sqrt{\ep}}
            -
                \hat{u}({ {1}}; \, {\blds \theta}) \,
                e^{-(1-x)/\sqrt{\ep}}
            , 
\end{equation}
where $\hat{u}$ is the network output as defined in (\ref{e:NN}).

By the construction (\ref{e:app_sol_enriched_rea-def}), the sl-FVNN prediction $\widetilde{u}$ 
satisfies the zero boundary value at $x = 0,1$, 
up to an exponentially {small  term}. 
{
Indeed, evaluating (\ref{e:app_sol_enriched_rea-def}) at the boundary points gives 
$$
\widetilde u(0;\theta) = -\widehat u(1;\theta)e^{-1/\sqrt{\varepsilon}}, \qquad \widetilde u(1;\theta) = -\widehat u(0;\theta)e^{-1/\sqrt{\varepsilon}}. 
$$
Thus, provided that the network output remains bounded, the boundary discrepancy is of order \(O(e^{-1/\sqrt{\varepsilon}})\), which is exponentially small as \(\varepsilon\to0\). 
Therefore, the homogeneous boundary conditions are enforced up to computationally negligible exponentially small errors.
}

{We } define  the sl-FVNN loss for this problem in the form, 
\begin{equation}\label{e:loss_PINNs_mod_rea-dif}
        \tilde{\mathcal{L}}_{FV}(\blds{\tht}; \mathcal{T}) 
        =
        \frac{1}{M-2} \sum_{i=2}^{M-1} 
        \Big| 
        - \ep d_x {\widetilde{u}} (x_{i+1}) 
        + \ep d_x {\widetilde{u}} (x_{i-1})
        + {\widetilde{u}} (x_i)  |K_i| 
        - f(x_i) |K_i|  
        \Big|^{2}. 
\end{equation}
\smallskip

We notice from the computational results in Section \ref{S.comp_RD_1D}   
below that 
our sl-FVNN outperforms the FVNN (as well as PINN and sl-PINN) especially when the parameter $\ep$ is small.

\subsection{Nonlinear Convection Diffusion Equation (CD eqn)}\label{S.cd_1d} 
With $g(u^\ep) = - d_x u^\ep  + (u^\ep)^3$ and $u_0= u_1 = 0$ in (\ref{e:1d_gen}), 
we consider the singularly perturbed convection diffusion equation with a non-linear reaction term:
\begin{equation}\label{e:cd_1d}
    \left\{
    \begin{array}{rl}
        \spacer
        - \ep d_x^2 u^\ep  - d_x u^\ep  + (u^\ep)^3   = f(x), & 0 < x <1, \\
        u^\ep = 0, & x=0, 1;
    \end{array}
    \right.
\end{equation}
see, e.g., \cite{JP07} for the vanishing diffusivity limit and asymptotic analysis for (\ref{e:cd_1d}). 
\smallskip

To derive a proper FVNN loss function for (\ref{e:cd_1d}), 
we integrate the equation over $K_i$, integrate it by parts, and find that 
\begin{equation}\label{e:1d_cd_w}
        - \ep d_x u^\ep(x_{i+1}) 
        + \ep d_x u^\ep(x_{i-1})       
        - u^\ep(x_{i+1}) 
        + u^\ep(x_{i-1})
        + \int_{K_i} (u^\ep)^3 \, dx
    =        
        \int_{K_i} f \, dx.
\end{equation}
Hence, using this proper treatment of the lower order term, $g(u^\ep) = - d_x u^\ep  + (u^\ep)^3$, we define the FVNN prediction $\hat{u}$ of (\ref{e:cd_1d}) by minimizing the loss, 
\begin{equation}\label{e:loss_FVNN_cd_1d}
    \begin{split}
        \mathcal{L}_{FV}(\blds{\tht}; \mathcal{T}) 
        =
        &
        \frac{1}{M-2} \sum_{i=2}^{M-1} 
        \Big| 
        - \ep d_x {\hat{u}} (x_{i+1}) 
        + \ep d_x {\hat{u}} (x_{i-1})
        { - {\hat{u}} (x_{i+1}) 
        + {\hat{u}} (x_{i-1})
        }\\
        &
        \qquad \qquad \qquad
        + \hat{u}^3 (x_i)  |K_i| 
        - f(x_i) |K_i|  
        \Big|^{2} \\
        &
        +
        | \hat{u} (0)  |^2  +
        | \hat{u} (1)  |^2.
    \end{split}
\end{equation}
\smallskip

Now, in order to construct the singular layer FVNN (sl-FVNN) for (\ref{e:cd_1d}),   
we first observe that
the limit problem at the vanishing 
diffusivity ($\ep = 0$) is given by 
\begin{equation}\label{e:rea_con_dif_eqn_nonlinear_limit}
    \left\{
    \begin{array}{rl}
        \spacer
        - d_x u^0  + (u^0)^3   = f(x), & 0 < x <1, \\
        u^0 = 0, & x= 1.
    \end{array}
    \right.
\end{equation} 

As verified in \cite{JP07}, the boundary layer associated with (\ref{e:cd_1d}) 
is linear although it is a non-linear problem.
In fact, by the matching asymptotic analysis,  
we see that 
the boundary layer of size $\ep$ occurs near the {\it outflow} boundary at $x = 0$.  
Further we find that 
the {asymptotic equation} for the corrector 
$\varphi \simeq u^\ep - u^0$  
is written  
in the form,
\begin{equation}\label{e:cor_eqn_con-rea-dif}
    \left\{
        \begin{array}{rl}
            \spacer
            -\ep d^2_x \varphi  - d_x \varphi = 0, 
            &
            0 < x < 1, \\
            \spacer
            \varphi= - u^0 , 
            &
            x = 0, \\
            {
            \varphi= 0 , }
            &
            {
            x = 1. }
        \end{array}
    \right.
\end{equation}
Solving the equation above, we obtain the explicit expression of $\varphi$ as
\begin{equation}\label{e:cor__con-rea-dif}
    \varphi(x)
        {
        =
            -u^0(0) \frac{e^{-x/\varepsilon}-e^{-1/\varepsilon}} {1-e^{-1/\varepsilon}}
        }
        =
            -u^0(0) \, e^{-x/\ep}
            + e.s.t.,  
\end{equation}
{
where \({e.s.t.}\) denotes terms that are exponentially small as \(\varepsilon\to0\).
}
The asymptotic analysis in \cite{JP07} confirms that   
the corrector $\varphi$ 
exhibits the singular behavior of $u^\ep$ at a small diffusivity $\ep$.

Using the profile of $\varphi$ in (\ref{e:cor__con-rea-dif}),  
we propose  
the sl-FVNN predicted solution $\tilde{u}$ for (\ref{e:cd_1d}) in the form,
\begin{equation}\label{e:app_sol_enriched_con-rea-def}
    \widetilde{u}(x ; \, {\blds \theta}) 
        =
            { (x-1)} \, 
            ( 
                \hat{u}({ {x}}; \, {\blds \theta})
            -
                \hat{u}({ {0}}; \, {\blds \theta}) \,
                e^{-x/\ep}
            ), 
\end{equation}
with the network output $\hat{u}$ as in (\ref{e:NN}).

The sl-FVNN prediction $\widetilde{u}$ in (\ref{e:app_sol_enriched_con-rea-def}) 
satisfies the zero boundary value at $x = 0,1$, 
up to an exponentially small error. 
Thus, we train the sl-FVNN predicted solution $\widetilde{u}$ by
minimizing the loss defined by 
\begin{equation}\label{e:loss_sl_FVNN_cd_1d}
    \begin{split}
        \tilde{\mathcal{L}}_{FV}(\blds{\tht}; \mathcal{T}) 
        =
        &
        \frac{1}{M-2} \sum_{i=2}^{M-1} 
        \Big| 
        - \ep d_x {\widetilde{u}} (x_{i+1}) 
        + \ep d_x {\widetilde{u}} (x_{i-1})
        { - {\widetilde{u}} (x_{i+1}) 
        + {\widetilde{u}} (x_{i-1})
        }\\
        &
        \qquad \qquad \qquad
        + {\widetilde{u}^3} (x_i)  |K_i| 
        - f(x_i) |K_i|  
        \Big|^{2}  
        .
    \end{split}
\end{equation}
\smallskip

We see below the computational results in Section \ref{S.comp_CD_1D}   
below that 
our sl-FVNN outperforms the FVNN (as well as PINN and sl-PINN) especially when the parameter $\ep$ is small.

\subsection{Stationary Burgers Equation (Burgers eqn)}\label{S.Burgers_1d}
In this section,  
we set 
$g(u^\ep) = u^\ep \, d_x u^\ep$ in (\ref{e:1d_gen}), $u_0 = -2$, and $u_1 = -1$,  
and consider the stationary 1D Burgers  equation:
\begin{equation}\label{e:Burgers}
\left\{\begin{array}{rl}
	\spacer
	-\ep d_x^2 u^\ep  + u^\ep \, d_x u^\ep  & = f(x), 
		 \quad 0 < x < 1,\\
         u^\ep & = {- 2}, \quad  x=0,  
        \\
        u^\ep & = - 1,  \quad x=1
        .
\end{array}\right.
\end{equation}
Here  
$0 < \epsilon \ll 1$ is a small viscosity parameter and 
$f$ is a smooth data. 
We have chosen the boundary values of $u^\ep$ at $x=0, 1$ 
as negative numbers $u_0 < u_1 < 0$.   
Under this setting,   
the convection $u^\ep u^\ep_x $ occurs always in one direction from right to left so that 
$u^\ep \leq -1$ in $(0, 1)$ 
and 
the boundary layer occurs near only   
the {\it outflow} boundary at $x=0$; 
see, e.g., \cite{CJL19, GJL22} for more information and detailed boundary layer analysis. 
\smallskip

To define the FVNN prediction $\hat{u}$,  
we integrate (\ref{e:Burgers}) over $K_i$, integrate it by parts, and write 
\begin{equation}\label{e:1d_Burgers_w}
        - \ep d_x u^\ep(x_{i+1}) 
        + \ep d_x u^\ep(x_{i-1})       
        + \dfrac{1}{2} (u^\ep)^2(x_{i+1}) 
        - \dfrac{1}{2} (u^\ep)^2(x_{i-1}) 
    =        
        \int_{K_i} f \, dx.
\end{equation}
Using this proper treatment of local fluxes, we define the FVNN prediction $\hat{u}$ of (\ref{e:1d_Burgers_w}) by minimizing the loss below:
\begin{equation}\label{e:loss_FVNN_Burgers}
    \begin{split}
        \mathcal{L}_{FV}(\blds{\tht}; \mathcal{T}) 
        =
        &
        \frac{1}{M-2} \sum_{i=2}^{M-1} 
        \Big| 
        - \ep d_x {\hat{u}} (x_{i+1}) 
        + \ep d_x {\hat{u}} (x_{i-1}) 
        + { \dfrac{1}{2} \hat{u}^2 (x_{i+1}) 
        - \dfrac{1}{2} \hat{u}^2 (x_{i-1})
        } \\
        &
        \qquad \qquad \qquad  
        - f(x_i) |K_i|  
        \Big|^{2} \\
        &
        +
        | \hat{u} (0) + 2 |^2  +
        | \hat{u} (1) + 1 |^2.
    \end{split}
\end{equation}
\smallskip

Now, for the singular layer FVNN (sl-FVNN) for (\ref{e:Burgers}),   
we notice that
the limit problem at the vanishing 
viscosity ($\ep = 0$) is  obtained 
by setting $\ep = 0$ in (\ref{e:Burgers}) 
and imposing the {\it inflow} boundary condition at $x=1$:
\begin{equation}\label{e:Burgers_limit}
	\left\{
	\begin{array}{l}
		\spacer
			u^0 u^0_x = f(x), \quad 0 < x < 1, \\
			u^0(1) = - 1.
	\end{array}
	\right.
\end{equation}
The formal limit $u^0$, a solution of (\ref{e:Burgers_limit}), is then found as   
\begin{equation}\label{e:u^0}
    u^0 (x) 
        =
            - \Big(
                2 \int_1^x f(x^\pri) \, dx^\pri 
                + 1
            \Big)^{\frac{1}{2}}. 
\end{equation}
\smallskip

Performing the matching asymptotics, 
we find that the size of boundary layer near the {\it outflow} boundary at $x = 0$ is 
of order $\ep$, and 
that  the asymptotic equation for   
$\varphi \simeq u^\ep - u^0$ is given by 
\begin{equation}\label{e:cor_eqn_Burgers}
    \left\{
        \begin{array}{rl}
            \spacer
            -\ep d^2_x \varphi  
            + 
            \big(
            u^0(0) + \varphi 
            \big)
            \pa_x \varphi   = 0, 
            &
            0 < x < 1, \\
            \spacer
            \varphi= - 2 - u^0 , 
            &
            x = 0, \\
            \varphi \rightarrow 0
            &
            \text{as }
            x \rightarrow 1.
        \end{array}
    \right.
\end{equation}
As shown in, e.g., \cite{CJL19, GJL22}, 
the corrector $\varphi$ is explicitly written as 
\begin{equation}\label{e:cor_Burgers}
    \varphi(x)
        =
            \dfrac{
                2u^0(0) \big(2 + u^0(0)\big)
                e^{u^0(0)x/{\ep}}
            }
            {
                2 - u^0(0) 
                -
                \big(2 + u^0(0)\big)
                e^{u^0(0)x/{\ep}}
            }
            .
\end{equation}
One can further verify that 
the fast decaying part of $u^\ep$ (near the {\it outflow} boundary at $x=0$) is 
described by the corrector $\varphi$, and 
that 
\begin{equation}\label{e:conv_bur}
    \|
        u^\ep - (u^0 + \varphi )
    \|_{L^2(0, 1)}
        \leq
            \kap \ep, 
    \qquad
    \|
        u^\ep - u^0  
    \|_{L^2(0, 1)}
        \leq
            \kap \ep^\frac{1}{2}, 
\end{equation}
for a constant $\kap >0$ independent of $\ep$. 
\smallskip

Using the network output $\hat{u}(x; \blds{\tht})$ as in (\ref{e:NN}), 
we set $\tilde{\varphi}(x; \blds{\tht})$ as the right-hand side of 
(\ref{e:cor_Burgers}) with $u^0(0)$ replaced by $\hat{u}(0; \blds{\tht})$, and  
introduce the normalization $\overline{\varphi}$ of $\tilde{\varphi}$ by {$\tilde{\varphi}(0; \blds{\tht})$}:  
\begin{equation}\label{e:cor_Burgers_prof}
    \overline{\varphi}(x; \blds{\tht}) = \dfrac{\tilde{\varphi}(x; \blds{\tht})}{\tilde{\varphi}(0; \blds{\tht})}, 
    \qquad
    \tilde{\varphi}(0; \blds{\tht})
        =
                {
            \dfrac{
                2 \hat{u}(0; \blds{\tht}) \big(2 + \hat{u}(0; \blds{\tht}) \big)
            }
            {
                2 - \hat{u}(0; \blds{\tht}) 
                -
                \big(2 + \hat{u}(0; \blds{\tht})\big)
                }
                }
            .
\end{equation}
We then define 
the sl-FVNN predicted solution $\tilde{u}$ for  (\ref{e:Burgers}) in the form,
\begin{equation}\label{e:app_sol_enriched_Burgers}
    \widetilde{u}(x ; \, {\blds \theta}) 
        =
            (x-1) \Big(
                        \hat{u}({ {x}}; \, {\blds \theta})
                        -
                        \hat{u}({ 0}; \, {\blds \theta}) \, 
                        \overline{\varphi} (x; {\blds \theta})  
                    \Big)
            - 2 + x.
\end{equation}
{
The form of (\ref{e:app_sol_enriched_Burgers}) is chosen to impose the inhomogeneous boundary
condition in (\ref{e:Burgers}) while embedding the boundary layer corrector. 
In fact, the term $-2+x$ 
is a simple lifting function satisfying 
$$
(-2+x)|_{x=0}=-2,\qquad (-2+x)|_{x=1}=-1, 
$$
and the factor \(x-1\) is a smooth cut-off function that vanishes at \(x=1\).
Other smooth lifting and cut-off functions with the same properties could be used instead,  because 
the resulting smooth difference can be learned by the network part $\hat{u}$.  
The main role of the solution structure (\ref{e:app_sol_enriched_Burgers}) is to embed the outflow boundary layer
corrector into the FVNN architecture while enforcing the prescribed boundary values.
Hence}  the sl-FVNN prediction $\widetilde{u}$ satisfies the boundary condition in (\ref{e:Burgers}), 
\begin{equation}\label{e:BC_sl_FVNN_Burgers}
    \widetilde{u}(0 ; \, {\blds \theta}) = -2, 
    \qquad
    \widetilde{u}(1 ; \, {\blds \theta}) = -1, 
\end{equation} 
{and we} train $\widetilde{u}$ by
using the loss defined by 
\begin{equation}\label{e:loss_sl_FVNN_Burgers}
    \begin{split}
        \tilde{\mathcal{L}}_{FV}(\blds{\tht}; \mathcal{T}) 
        =
        &
        \frac{1}{M-2} \sum_{i=2}^{M-1} 
        \Big| 
        - \ep d_x {\widetilde{u}} (x_{i+1}) 
        + \ep d_x {\widetilde{u}} (x_{i-1})
        + { \dfrac{1}{2} {\widetilde{u}^2} (x_{i+1}) 
        - \dfrac{1}{2} {\widetilde{u}^2} (x_{i-1})
        } \\
        &
        \qquad \qquad \qquad
        - f(x_i) |K_i|  
        \Big|^{2}   .
    \end{split}
\end{equation}
\smallskip

We notice from the numerical results in Section \ref{S.comp_Burgers}   
below that 
our sl-FVNN outperforms the FVNN (as well as PINN and sl-PINN) especially when the parameter $\ep$ is small.

\section{Numerical Results}\label{S.Comp_1d}

In this section, 
we present numerical experiments demonstrating that the sl-FVNNs achieve superior performance compared to the standard FVNNs, in particular for the small perturbation parameter regime ($0 < \varepsilon \ll 1$), while maintaining efficient computational cost. 
To enable a direct comparison with the PINNs and sl-PINNs studied in \cite{GHJ} as well,  
we adopt the same computational setup (comparable to that in \cite{GHJ}) for all the numerical experiments performed in this article:
\begin{enumerate}
\item Number of epochs: EPOCHS = 2000
\item Number of layers: $L$ = 3 (single hidden layer)
\item Number of neurons: $N^l$ = 20
\item Number of collocation points: 50 
\item Activation function: sigmoid
\item Optimizer: LBFGS in PyTorch
\item The maximum number of iterations is set to 50,000 (2000 epochs × 25 LBFGS iterations), with
automatic early termination triggered by a tolerance threshold of $10^{-8}$
\end{enumerate}    

We test the model equations below for various perturbation parameters $\ep = 10^{-k}$, $4 \leq k \leq 12$ using the same setup described above. 

Because of the random initialization of weights and biases, standard FVNNs may fail to converge when $\ep$ is small; see, e.g., \cite{LMMK21, GLK19}.  
Hence, we perform five independent runs for each example and report the average outcome as our final result. 

{The plots shown below display raw approximations without any normalization, and they are selected from the five independent runs.  
For the comparisons between FVNNs and sl-FVNNs in the tables, we use the normalized relative error defined as}
{ 
\begin{equation}\label{e:rel_error}
\text{ Relative $L^{p}$ error of $u$} 
    = %
        \frac{\| u_{exact} - u \|_{L^p}}{\| u_{exact} \|_{L^p}},
\quad 
p = 2, \infty,  
\end{equation}
where $u_{exact}$ denotes the exact solution. 

Once the training process is complete,  
we evaluate the accuracy using testing points sampled more inside of a boundary layer, i.e., 
$50$ points inside of a boundary layer (of size $\ep^{1/2}$ or $\ep$)
and 500 points elsewhere.

\subsection{Reaction Diffusion Equation (RD eqn.)}\label{S.comp_RD_1D}
\noindent\\
For the computations of (\ref{e:rd_1D_eqn}), 
we introduce an exact solution $u_{exact}$  in the form, 
\begin{equation}\label{exact_rd}
    u_{exact}(x) 
        = 
        {
        (1-2 \ep) \bigg(
        1 - \dfrac{1 - e^{-1/\sqrt{\ep}}}{1 - e^{-2/\sqrt{\ep}}} 
        \Big(
        e^{-x/\sqrt{\ep}} + e^{-(1-x)/\sqrt{\ep}}
        \Big)
        \bigg) + x(1-x),
        }
\end{equation}
with the data, 
\begin{equation}\label{exact_rd_data}
    f(x) 
        = {
        1 + x(1-x).
        }
\end{equation}

Figure \ref{f:loss_rd_ep8} illustrates the loss values of both FVNN and sl-FVNN during the training at $\ep = 10^{-8}$. 
The graphs of exact solutions, FVNN predictions, and sl-FVNN predictions at different $\ep$ values appear in Figures \ref{f:rd_FVNN} and \ref{f:rd_slFVNN}. 
The comparison of numerical errors and training times between FVNNs and sl-FVNNs are summarized in Tables \ref{table1_rd} and \ref{table2_rd} respectively.

The numerical experiments show that sl-FVNNs significantly outperform standard FVNNs as the perturbation parameter $\ep>0$ gets smaller. 
While FVNNs experience slow and unstable training as $\ep$ decreases, sl-FVNNs converge rapidly due to the embedded singular-layer corrector, which provides an accurate analytic information of the boundary layer structure. 

Solution and error plots in Figures \ref{f:rd_FVNN} and \ref{f:rd_slFVNN} confirm that FVNNs fail to resolve sharp boundary layers, producing large spikes in $L^\infty$ error for $\ep \leq 10^{-6}$. 
In contrast, sl-FVNNs maintain uniformly high accuracy across all tested $\ep$ values. Quantitatively, sl-FVNNs achieve error reductions of up to $10^{-2} - 10^{-4}$ compared with FVNNs, as shown in Table \ref{table1_rd}.

From Tables \ref{table1_rd} and \ref{table1_rd_PINN}, we observe that, at $\ep = 10^{-8}$, the sl-FVNN and sl-PINN methods yield comparable accuracy, both achieving relative $L^2$ errors on the order of $10^{-4}$, demonstrating their effectiveness in resolving sharp layers over the PINNs and FVNNs.

Training-time comparisons in Table \ref{table2_rd} further demonstrate the efficiency of the proposed method, with sl-FVNNs generally requiring fewer iterations and achieving faster convergence. Overall, the results verify that sl-FVNNs provide both superior accuracy and improved computational performance, particularly in the singularly perturbed regime.

\begin{figure}
\centering
\subfloat[FVNNs]{\includegraphics[width=0.5\linewidth]{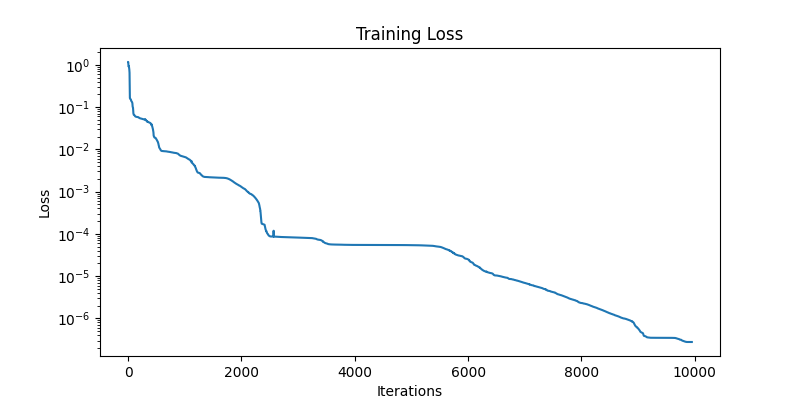}}
\hfill
\subfloat[sl-FVNNs]{\includegraphics[width=0.5\linewidth]{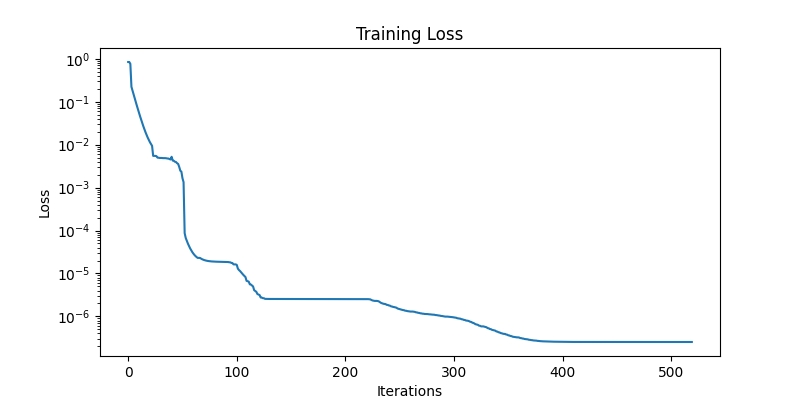}}
\caption{
{Loss values during the training of RD eqn. at $\ep = 10^{-8}$.
}
}\label{f:loss_rd_ep8}
\end{figure}

\begin{figure}
\centering
\subfloat[$\ep = 10^{-4}$]{\includegraphics[width=0.3\linewidth]{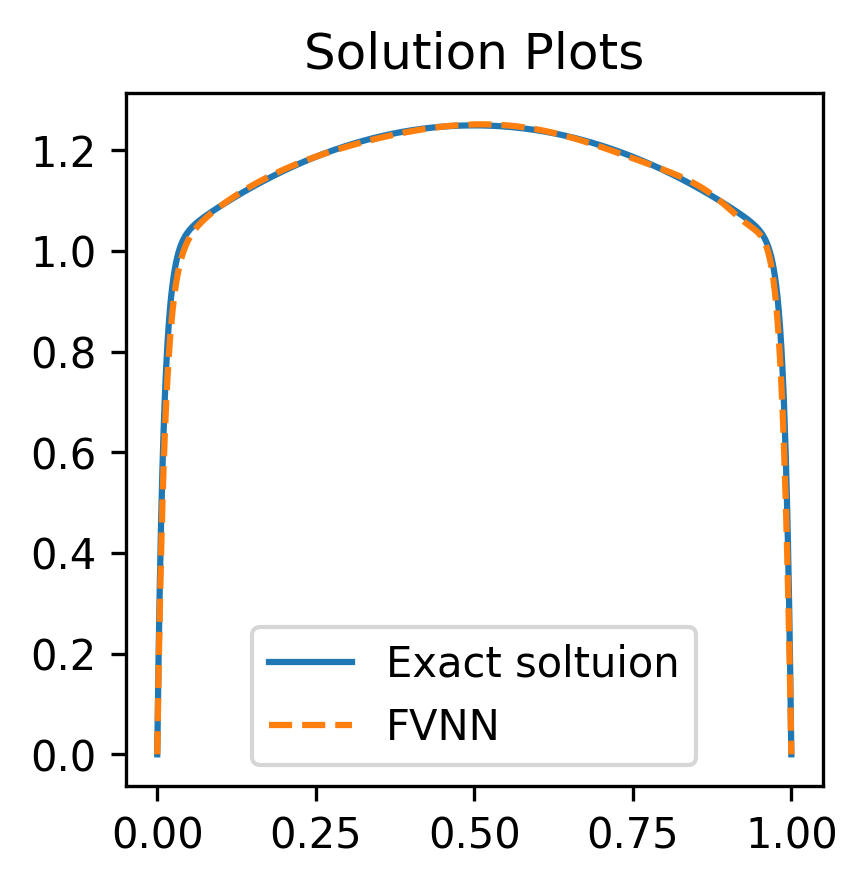}}
\subfloat[$\ep = 10^{-8}$]{\includegraphics[width=0.3\linewidth]{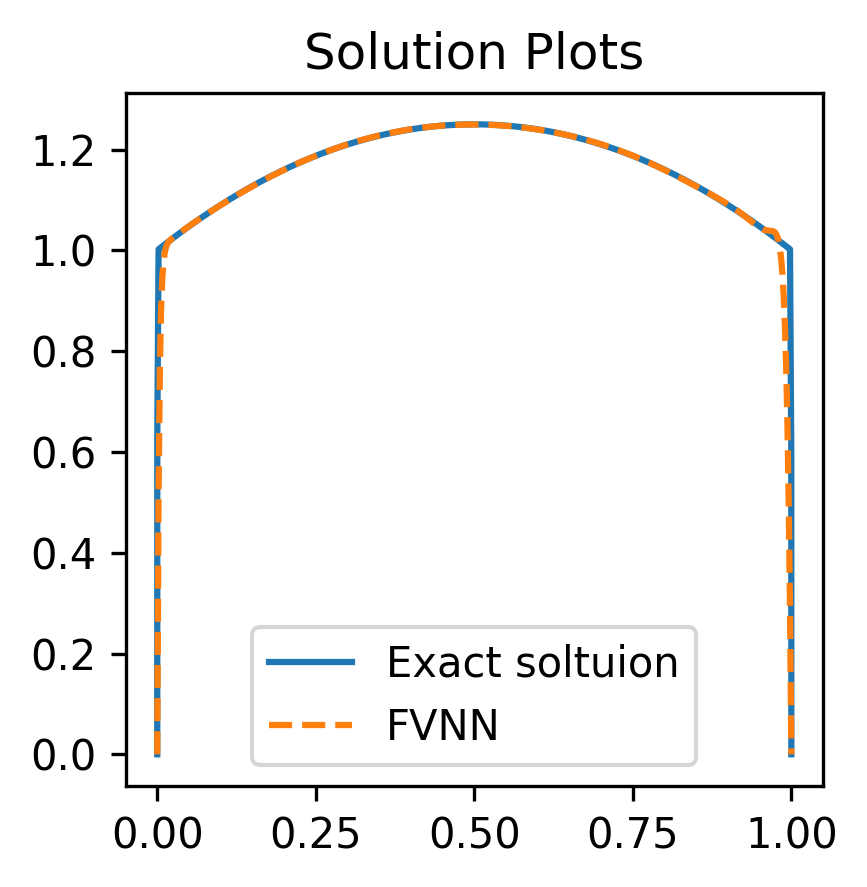}}
\subfloat[$\ep = 10^{-12}$]{\includegraphics[width=0.3\linewidth]{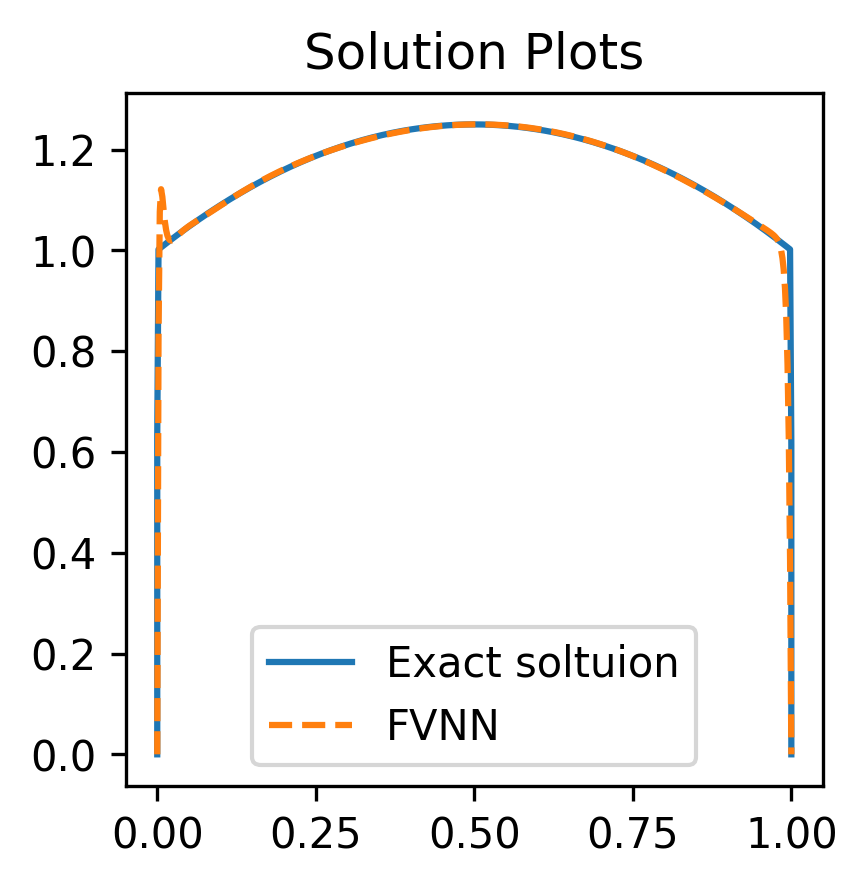}}
\vfill
\subfloat[$\ep = 10^{-4}$]{\includegraphics[width=0.3\linewidth]{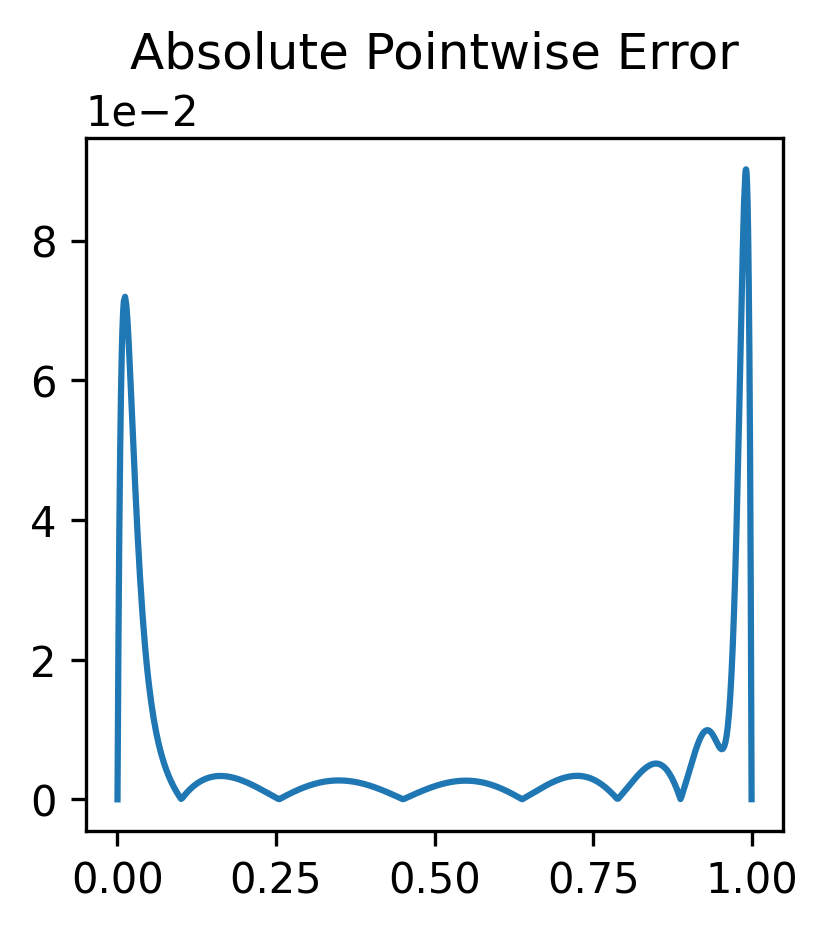}}
\subfloat[$\ep = 10^{-8}$]{\includegraphics[width=0.3\linewidth]{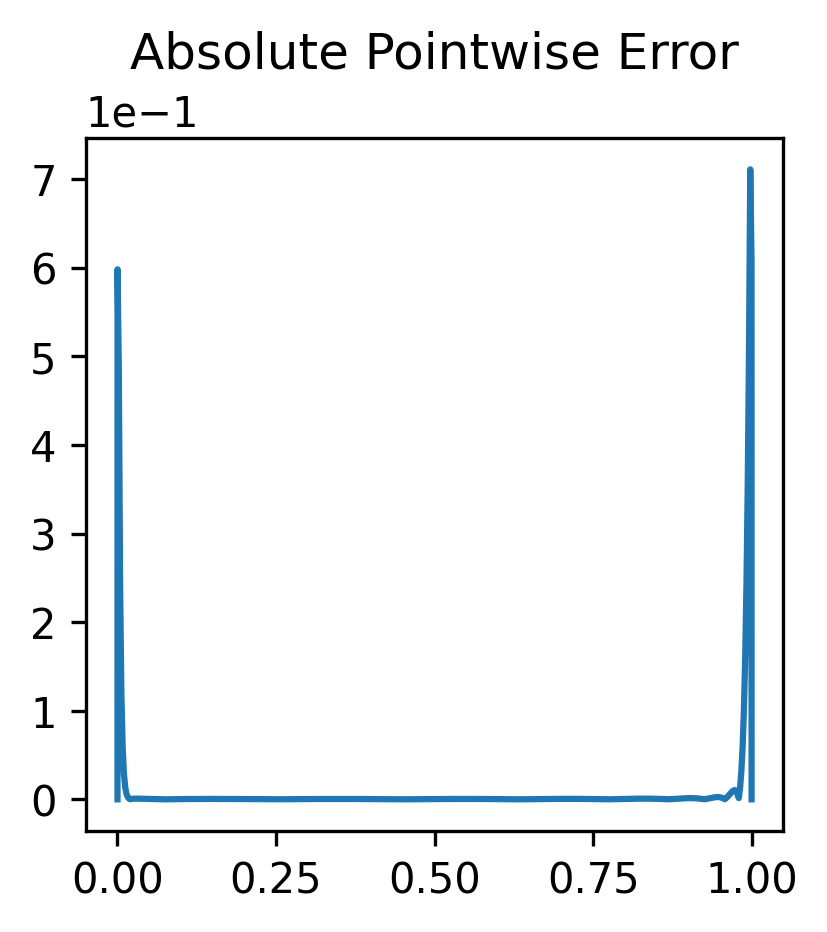}}
\subfloat[$\ep = 10^{-12}$]{\includegraphics[width=0.3\linewidth]{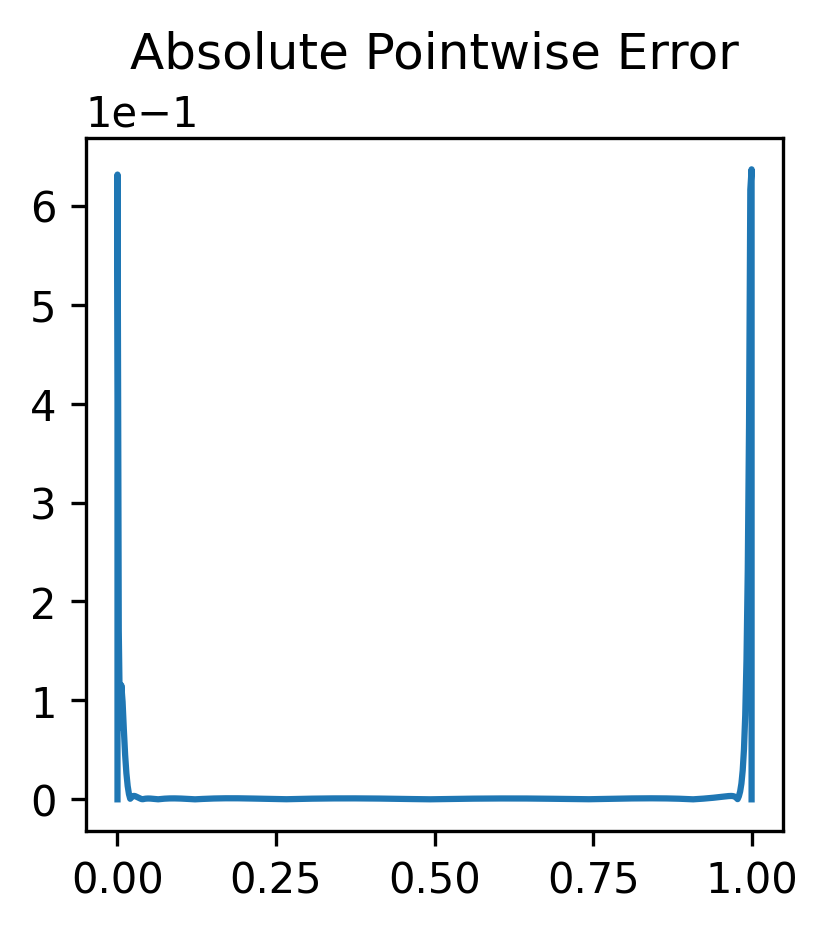}}
\caption{Exact solutions and FVNN predictions of RD eqn.}\label{f:rd_FVNN}
\end{figure}

\begin{figure}
\centering
\subfloat[$\ep = 10^{-4}$]{\includegraphics[width=0.3\linewidth]{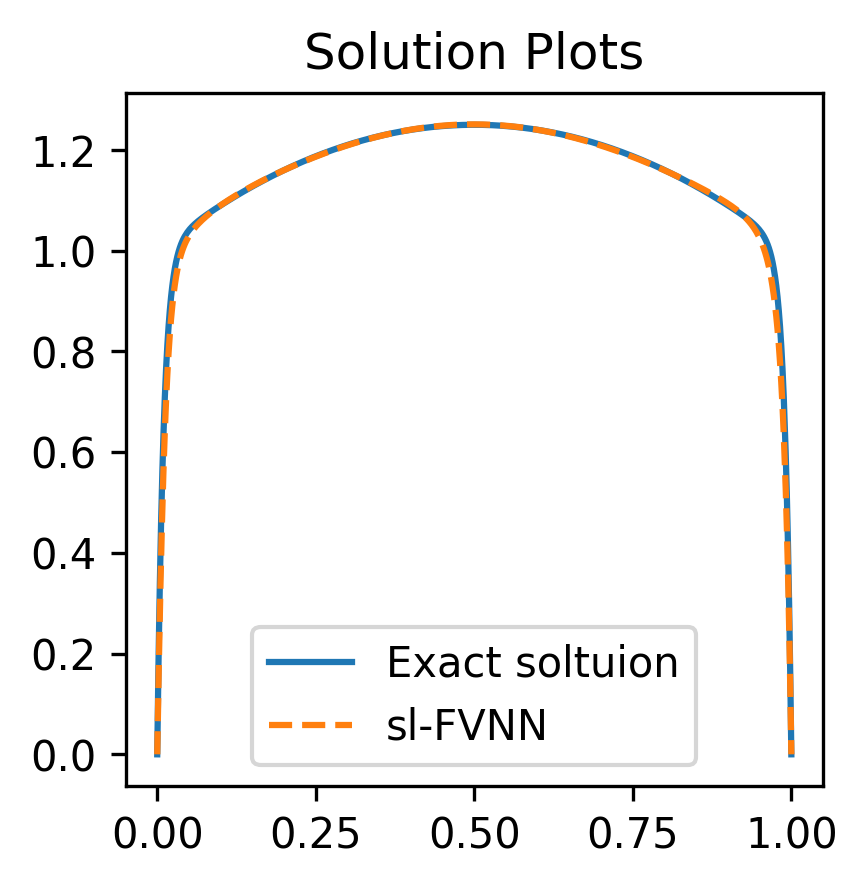}}
\subfloat[$\ep = 10^{-8}$]{\includegraphics[width=0.3\linewidth]{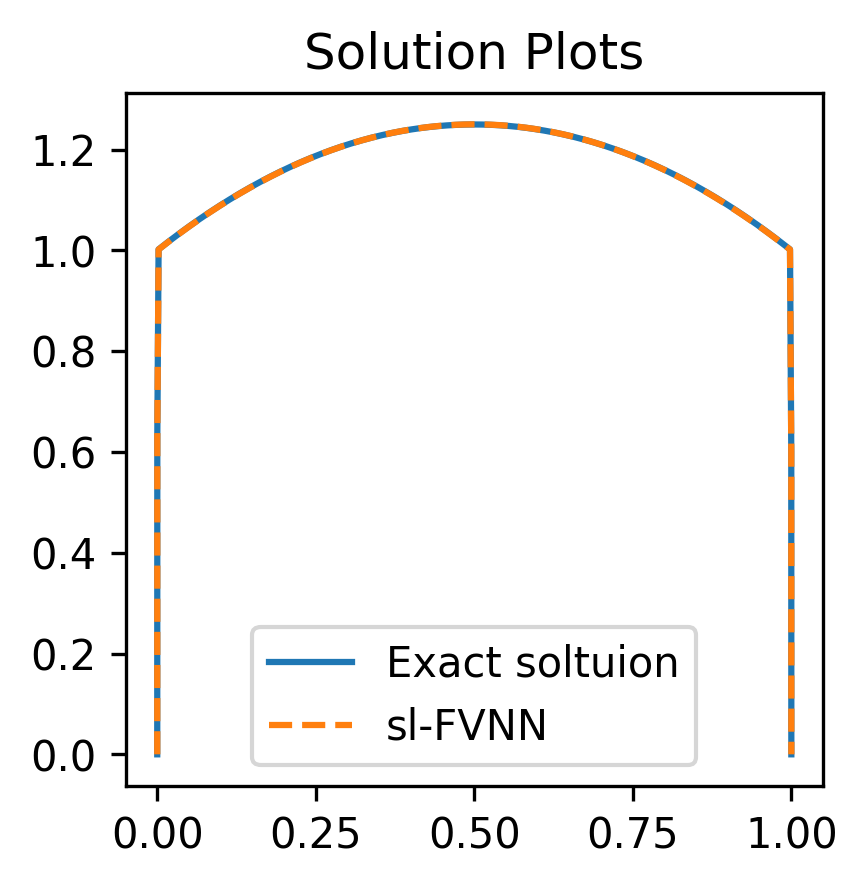}}
\subfloat[$\ep = 10^{-12}$]{\includegraphics[width=0.3\linewidth]{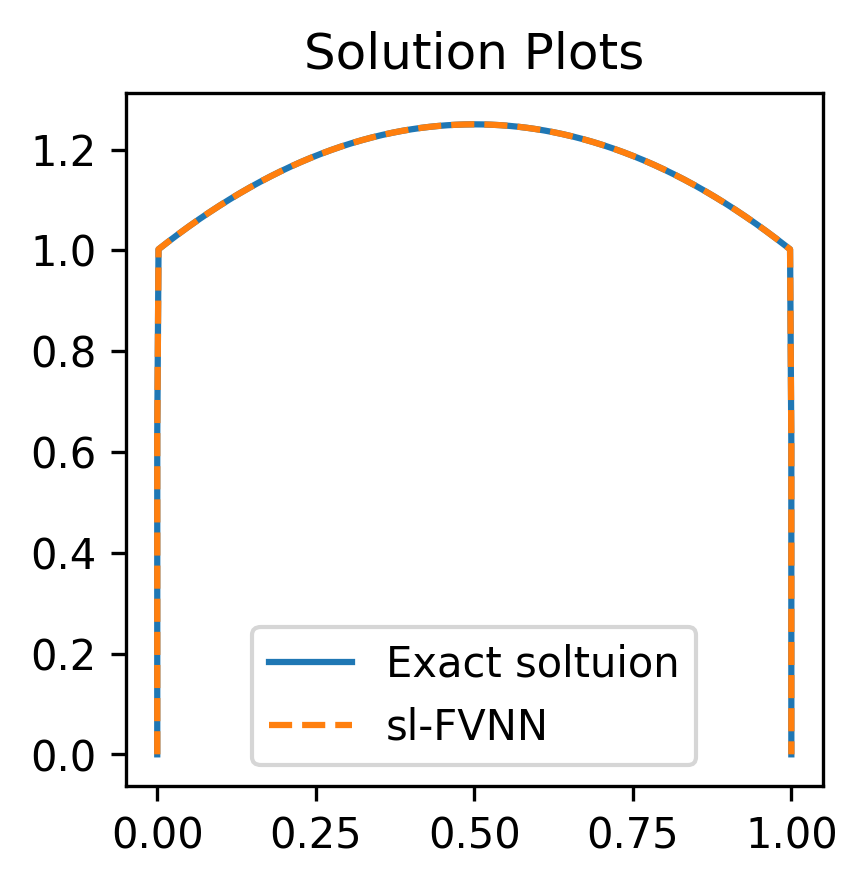}}
\vfill
\subfloat[$\ep = 10^{-4}$]{\includegraphics[width=0.3\linewidth]{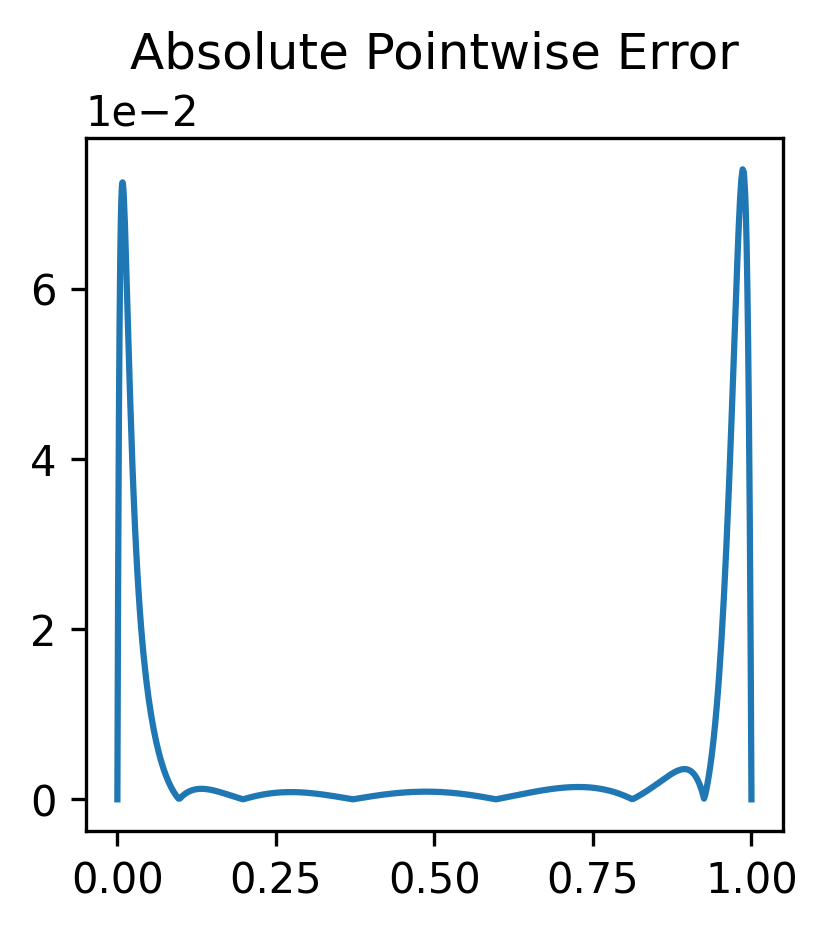}}
\subfloat[$\ep = 10^{-8}$]{\includegraphics[width=0.3\linewidth]{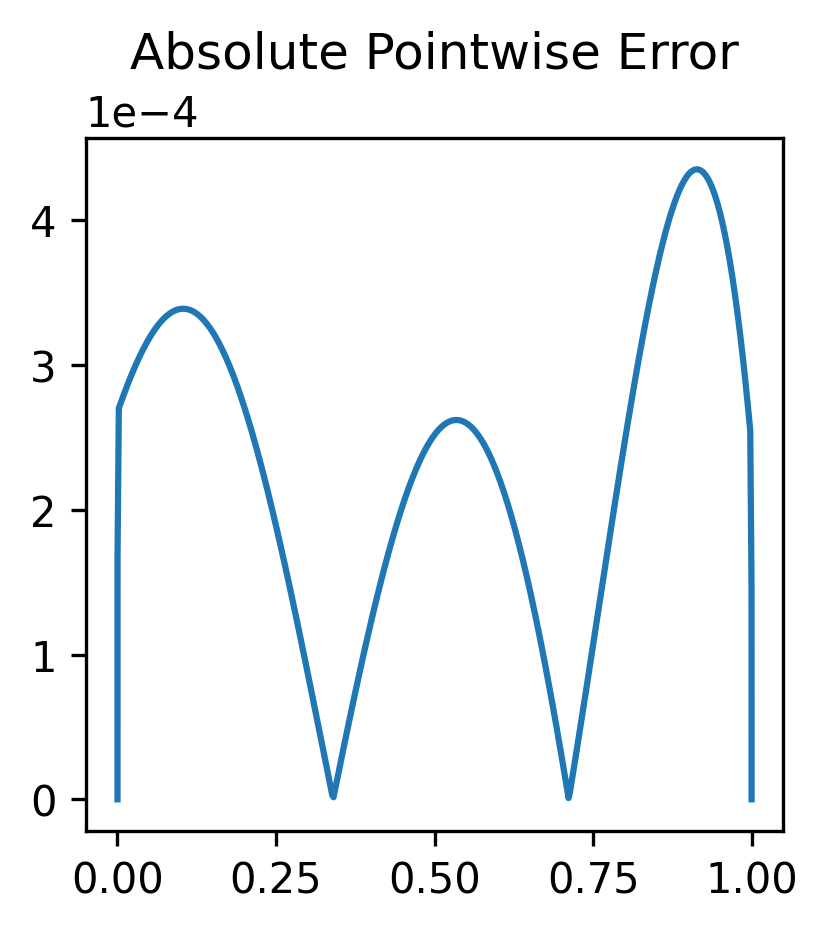}}
\subfloat[$\ep = 10^{-12}$]{\includegraphics[width=0.3\linewidth]{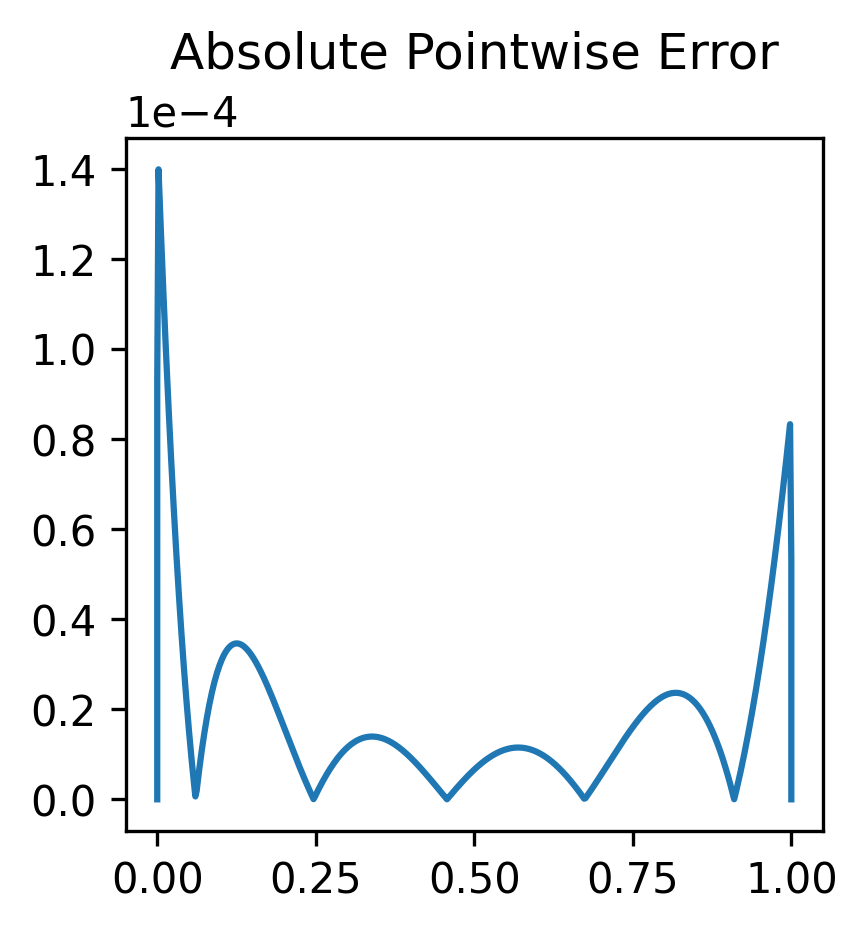}}
\caption{Exact solutions and sl-FVNN predictions of RD eqn.}\label{f:rd_slFVNN}
\end{figure}

\begin{table}[]
\begin{tabular}{|l|ll|ll|}
\hline
                      & FVNNs                                       &                             & {sl-FVNNs}                                    &                             \\ \hline
                      & \multicolumn{1}{l|}{Relative $L^{2}$ error} & Relative $L^{\infty}$ error & \multicolumn{1}{l|}{Relative $L^{2}$ error} & Relative $L^{\infty}$ error \\ \hline 
$\varepsilon=10^{-4}$ & 
\multicolumn{1}{l|}
{2.6209e-02}             & 7.2176e-02                  & \multicolumn{1}{l|}
{2.3589e-02}             & 5.9231e-02                  \\
$\varepsilon=10^{-6}$ & 
\multicolumn{1}{l|}
{1.1785e-01}              & 4.4217e-01                  & \multicolumn{1}{l|}
{3.9545e-03}             & 1.1379e-02                  \\
$\varepsilon=10^{-8}$ & 
\multicolumn{1}{l|}
{1.5613e-01}              & 5.6885e-01                  & \multicolumn{1}{l|}
{2.2174e-04}             & 3.4817e-04                  \\
$\varepsilon=10^{-10}$ & 
\multicolumn{1}{l|}
{1.6014e-01}              & 5.4955e-01                  & \multicolumn{1}{l|}
{3.6454e-05}             & 1.3175e-04                  \\
$\varepsilon=10^{-12}$ & 
\multicolumn{1}{l|}
{1.3276e-0}              & 5.0938e-01                  & \multicolumn{1}{l|}
{2.8209e-05}             & 1.1187e-04                  \\
\hline
\end{tabular}
\caption{Error comparison between FVNNs and {sl-FVNNs} for  RD eqn.}
\label{table1_rd}
\end{table}

\begin{table}[]
\begin{tabular}{|l|l|l|}
\hline
                      & 
                      PINNs (Relative $L^{2}$ error)
                      & 
                      {sl-PINNs} (Relative $L^{2}$ error)
                      \\ \hline
                      $\ep = 10^{-8}$
                      & 3.140e-01
                      & 1.245e-04 
                      \\ \hline  
\end{tabular}
\caption{
{
Error comparison between PINNs and {sl-PINNs} for  RD eqn. cited from \cite{GHJ}.
}
}
\label{table1_rd_PINN}
\end{table}

\begin{table}[]
\begin{tabular}{|l|ll|ll|}
\hline
                      & FVNNs                                       &                             & {sl-FVNNs}                                    &                             \\ \hline
                      & \multicolumn{1}{l|}{Training time} & Iterations/sec & \multicolumn{1}{l|}{Training time} & Iterations/sec \\ \hline 
$\varepsilon=10^{-4}$ & 
\multicolumn{1}{l|}
{9.02 sec}             & 229.28                  & \multicolumn{1}{l|}
{5.03 sec}             & 648.03                  \\
$\varepsilon=10^{-6}$ & 
\multicolumn{1}{l|}
{6.39 sec}              & 367.82                  & \multicolumn{1}{l|}
{1.14 sec}             & 422.44                  \\
$\varepsilon=10^{-8}$ & 
\multicolumn{1}{l|}
{2.49 sec}              & 3994.05                  & \multicolumn{1}{l|}
{2.62 sec}             & 198.36                  \\
$\varepsilon=10^{-10}$ & 
\multicolumn{1}{l|}
{2.85 sec}              & 2064.15                  & \multicolumn{1}{l|}
{1.17 sec}             & 604.80                  \\
$\varepsilon=10^{-12}$ & 
\multicolumn{1}{l|}
{6.78 sec}              & 1025.31                  & \multicolumn{1}{l|}
{1.43 sec}             & 1059.99                  \\
\hline
\end{tabular}
\caption{Training time comparison between FVNNs and {sl-FVNNs} for  RD eqn.}
\label{table2_rd}
\end{table}

\newpage

\subsection{Nonlinear Convection Diffusion equation (CD eqn.)}\label{S.comp_CD_1D}
\noindent\\
For the computations of (\ref{e:cd_1d}), 
we introduce an exact solution $u_{exact}$  in the form, 
\begin{equation}\label{exact_cd}
    u_{exact}(x) 
        = 
        (-1 + e^{-1/\ep})x^3 + 1 - e^{-x/\ep},  
\end{equation}
where the data $f$ is properly chosen by computing the right-hand side of (\ref{e:cd_1d}) with $u = u_{exact}$. 

The loss values of FVNN and sl-FVNN  at $\ep = 10^{-8}$ appear in 
Figure \ref{f:loss_cd_ep8} and   
the graphs of FVNN predictions and sl-FVNN predictions at different $\ep$ values appear in Figures \ref{f:cd_FVNN} and \ref{f:cd_slFVNN}. 
The comparison of numerical errors and training times between FVNNs and sl-FVNNs are summarized in Tables \ref{table1_cd} and \ref{table2_cd} respectively.

For the convection–diffusion equation (CD eqn), 
sl-FVNNs consistently outperform FVNNs in both accuracy and efficiency. 
FVNNs fail to resolve sharp boundary layers and exhibit unstable behavior as $\ep$ decreases, with relative errors growing by several orders of magnitude. 
In contrast, sl-FVNNs maintain uniform accuracy ($\simeq$ 1e-04) across all $\ep$ values, as shown in 
Figures \ref{f:cd_FVNN} - \ref{f:cd_slFVNN} and 
Table \ref{table1_cd}.

From Tables \ref{table1_cd} and \ref{table1_cd_PINN}, we observe that  
both singular layer enhanced methods (sl-PINNs and sl-FVNNs) significantly outperform their uncorrected counterparts (PINNs and FVNNs). 
Moreover, sl-FVNNs exhibit superior performance over sl-PINNs, achieving a relative 
$L^2$ error that is smaller by a factor of $10^2$.

Training curves (Figure \ref{f:loss_cd_ep8}) and 
Table \ref{table2_cd} confirm that sl-FVNNs also converge significantly faster, requiring less than one second of training time in most cases. These results highlight that sl-FVNNs offer a robust and scalable solution for handling stiff convection-dominated problems with steep boundary layers.

\begin{figure}
\centering
\subfloat[FVNNs]{\includegraphics[width=0.5\linewidth]{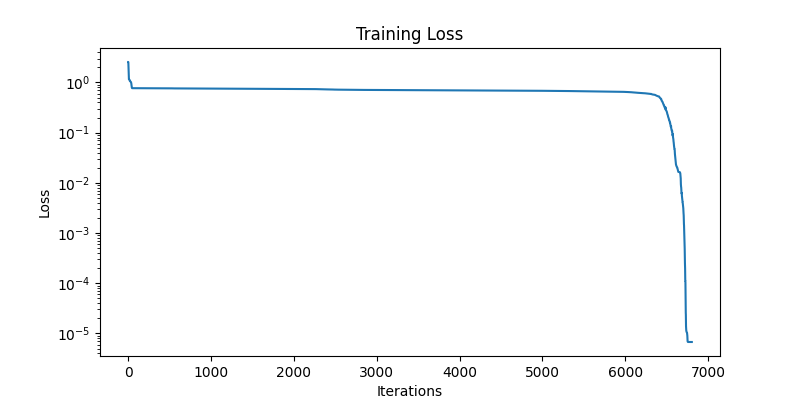}}
\hfill
\subfloat[sl-FVNNs]{\includegraphics[width=0.5\linewidth]{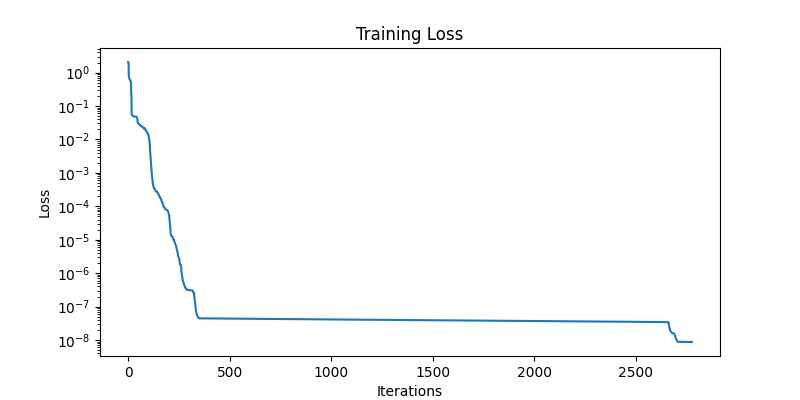}}
\caption{
{
Loss values during the training of CD eqn. at $\ep = 10^{-8}$.
}
}\label{f:loss_cd_ep8}
\end{figure}

\begin{figure}
\centering
\subfloat[$\ep = 10^{-4}$]{\includegraphics[width=0.3\linewidth]{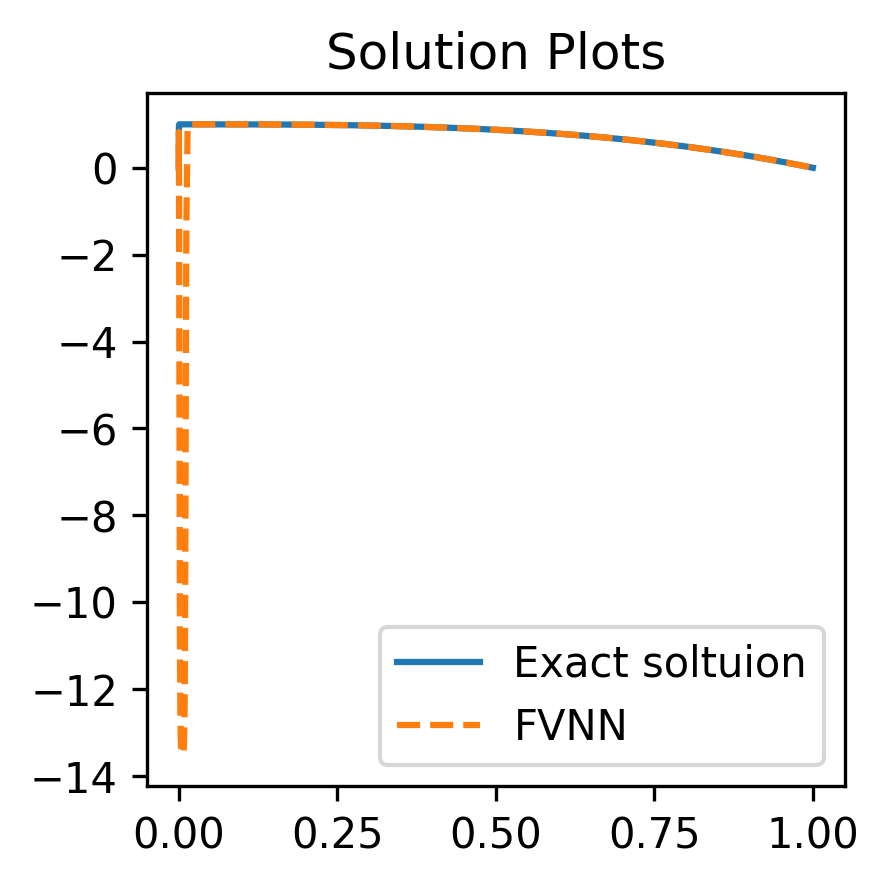}}
\subfloat[$\ep = 10^{-8}$]{\includegraphics[width=0.3\linewidth]{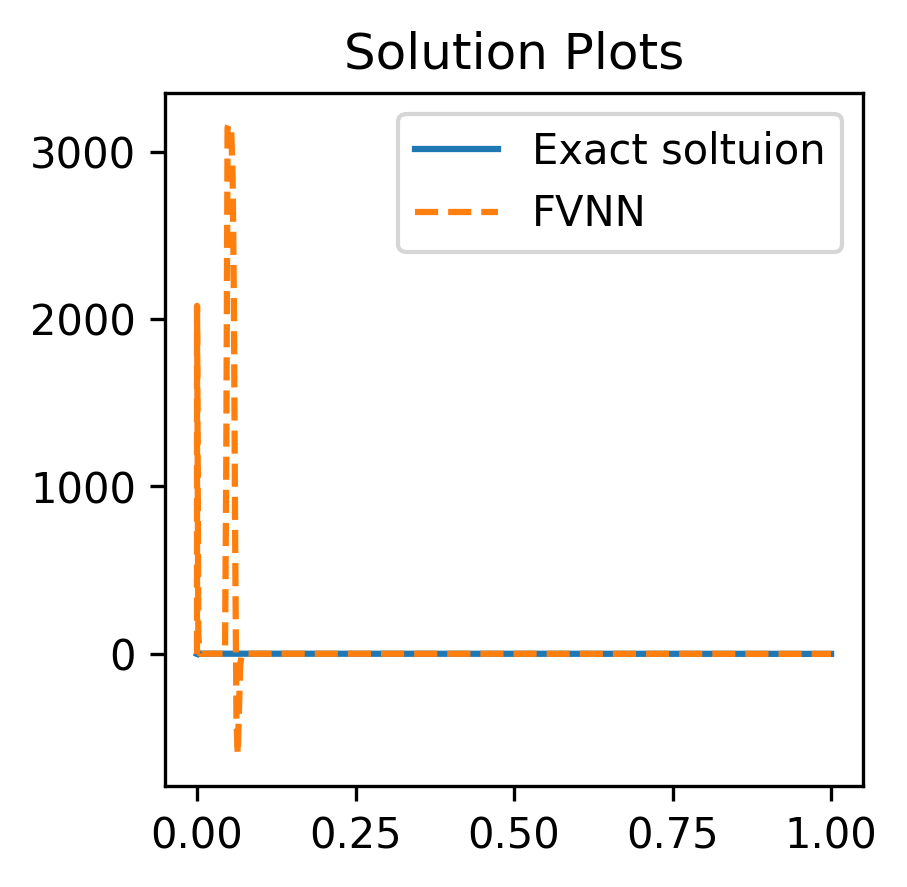}}
\subfloat[$\ep = 10^{-12}$]{\includegraphics[width=0.3\linewidth]{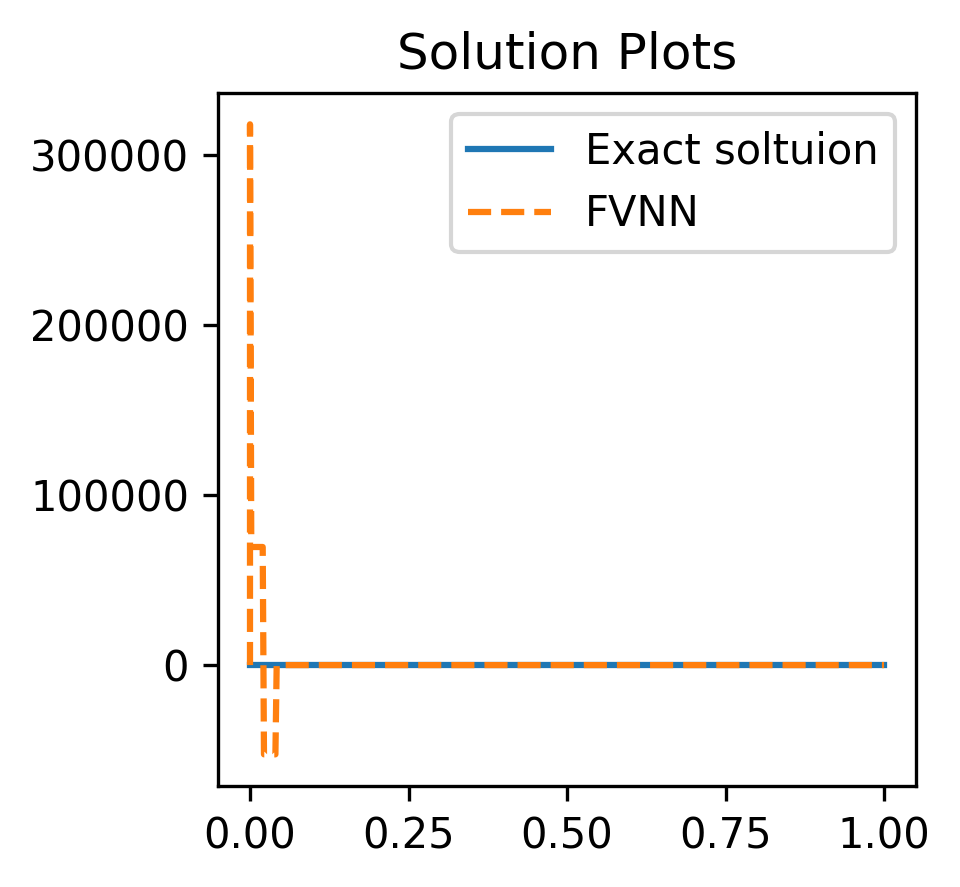}}
\vfill
\subfloat[$\ep = 10^{-4}$]{\includegraphics[width=0.3\linewidth]{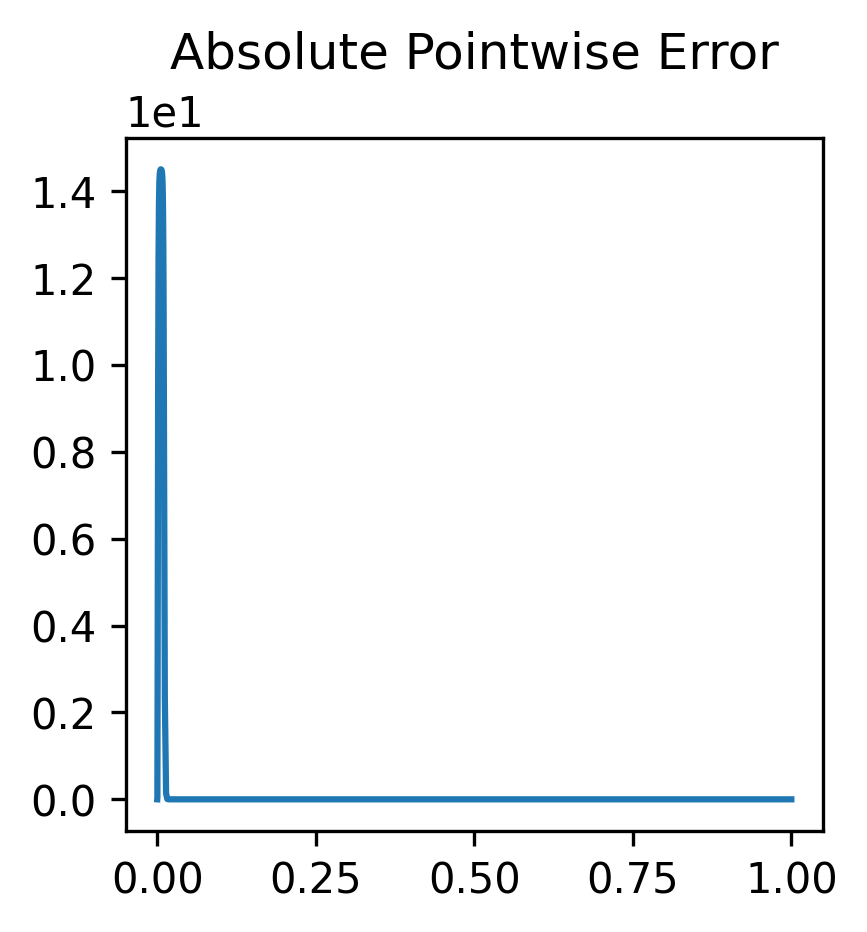}}
\subfloat[$\ep = 10^{-8}$]{\includegraphics[width=0.3\linewidth]{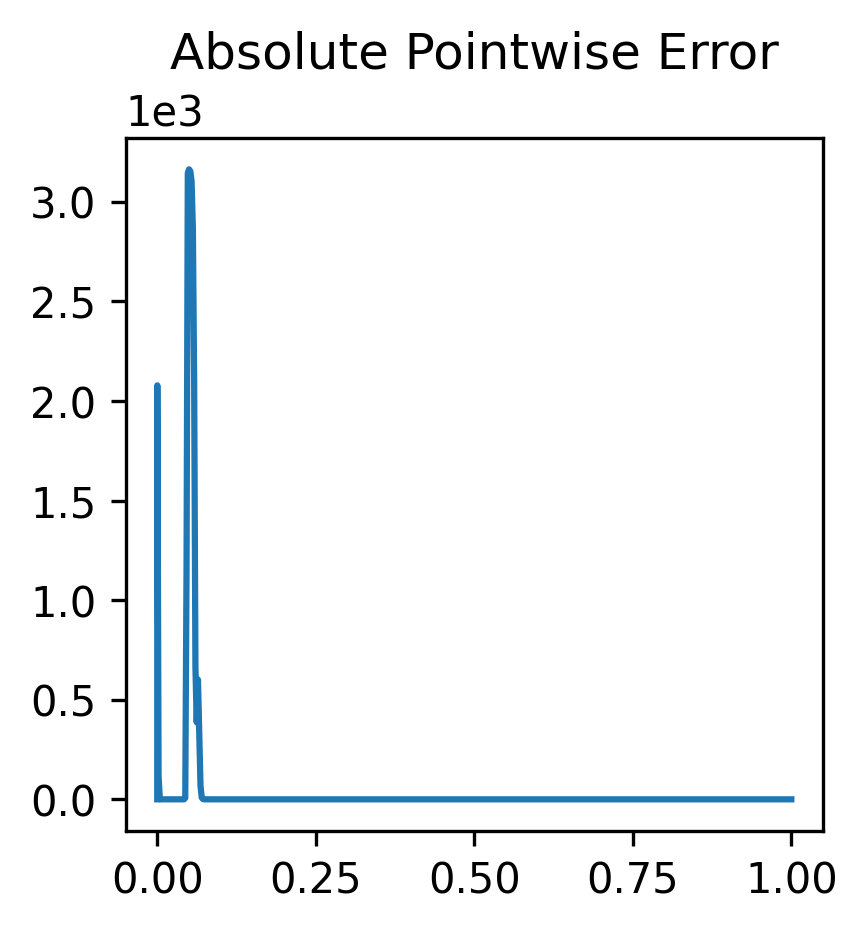}}
\subfloat[$\ep = 10^{-12}$]{\includegraphics[width=0.3\linewidth]{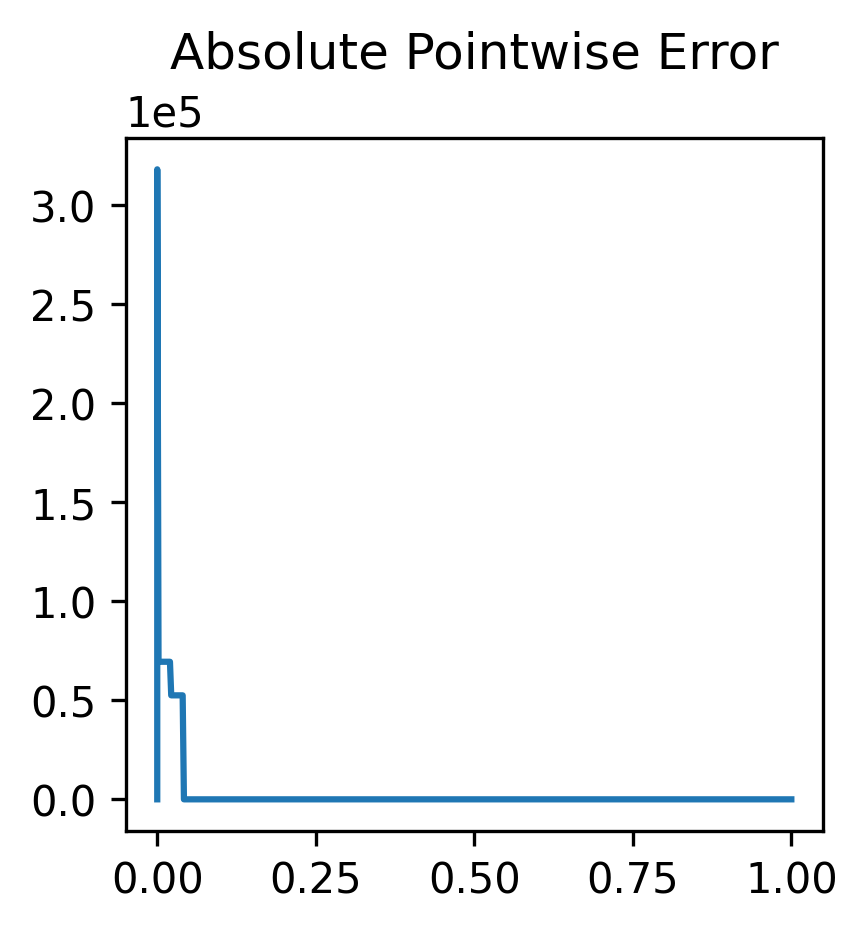}}
\caption{Exact solutions and FVNN predictions of CD eqn.}\label{f:cd_FVNN}
\end{figure}

\begin{figure}
\centering
\subfloat[$\ep = 10^{-4}$]{\includegraphics[width=0.3\linewidth]{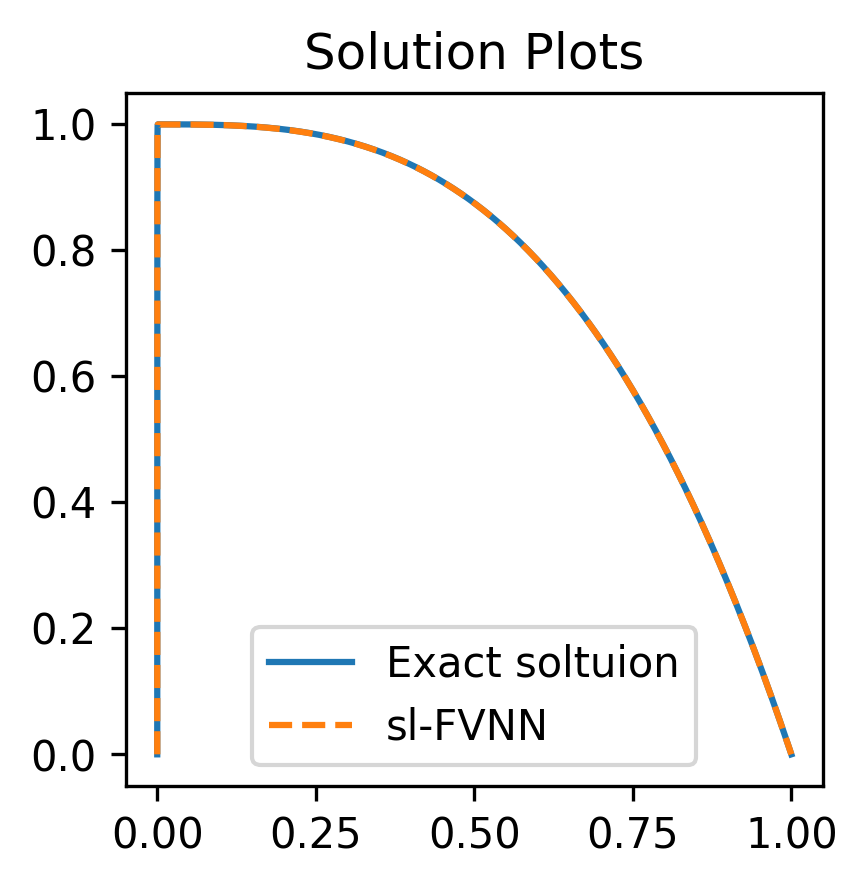}}
\subfloat[$\ep = 10^{-8}$]{\includegraphics[width=0.3\linewidth]{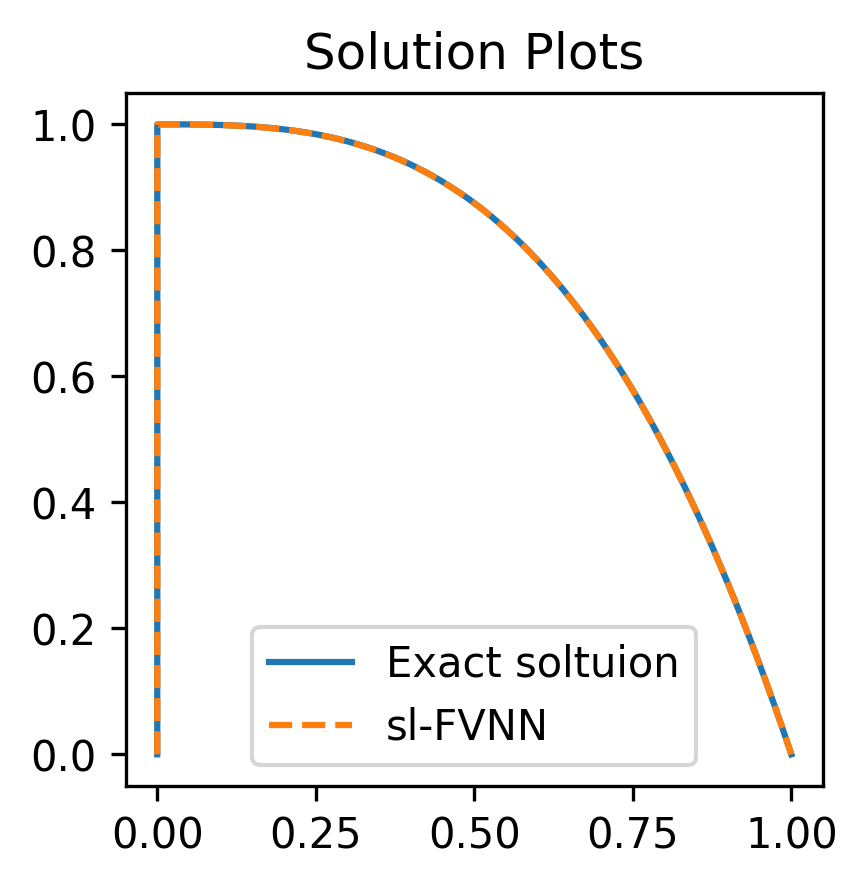}}
\subfloat[$\ep = 10^{-12}$]{\includegraphics[width=0.3\linewidth]{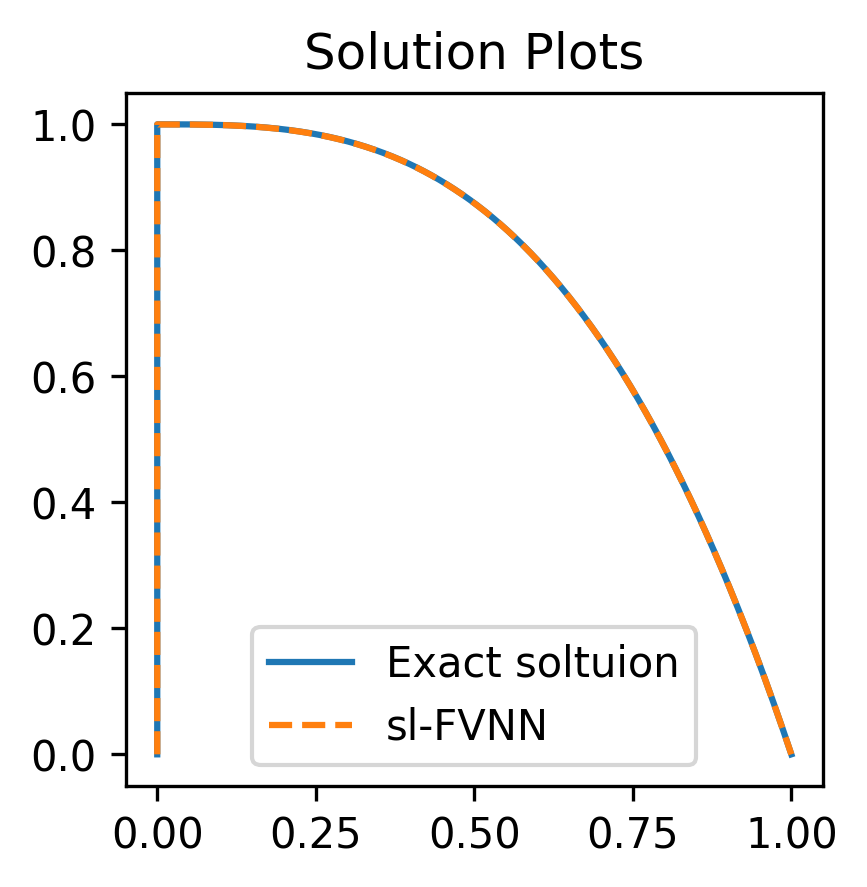}}
\vfill
\subfloat[$\ep = 10^{-4}$]{\includegraphics[width=0.3\linewidth]{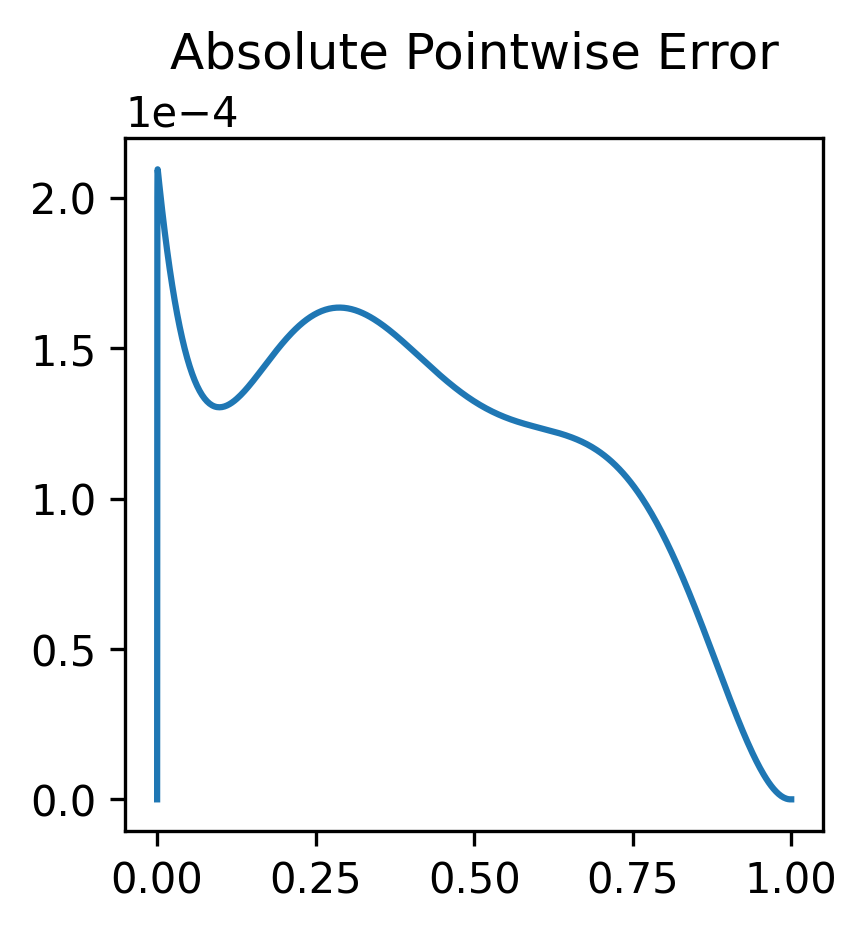}}
\subfloat[$\ep = 10^{-8}$]{\includegraphics[width=0.3\linewidth]{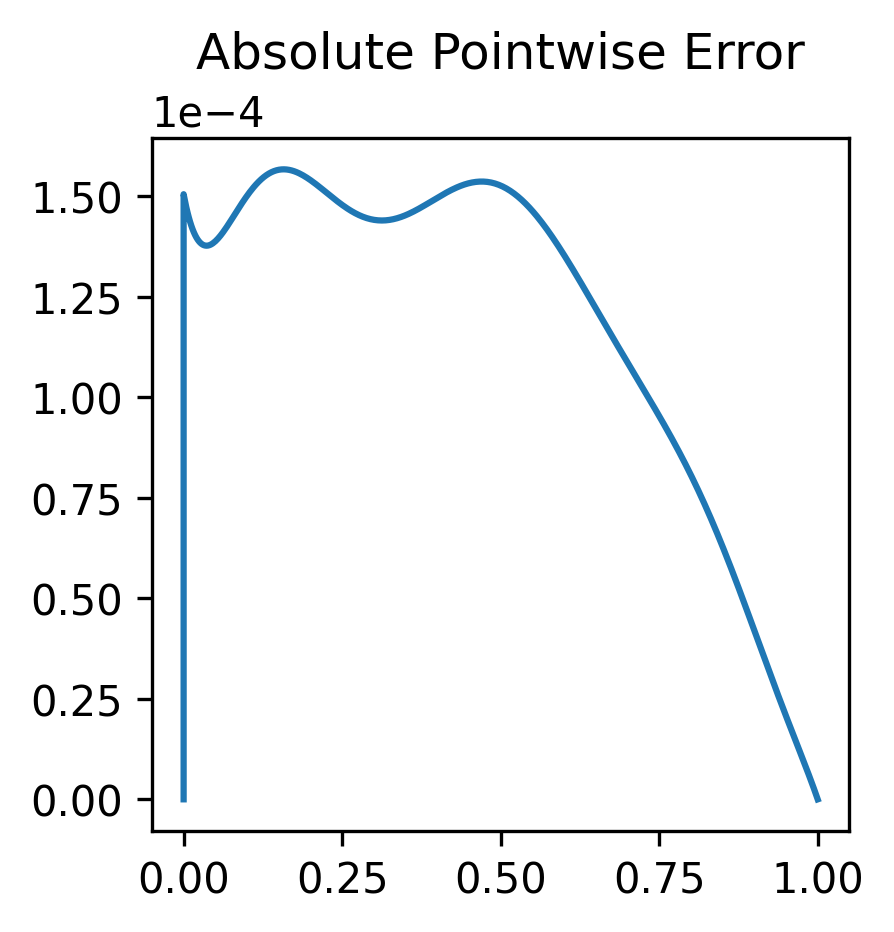}}
\subfloat[$\ep = 10^{-12}$]{\includegraphics[width=0.3\linewidth]{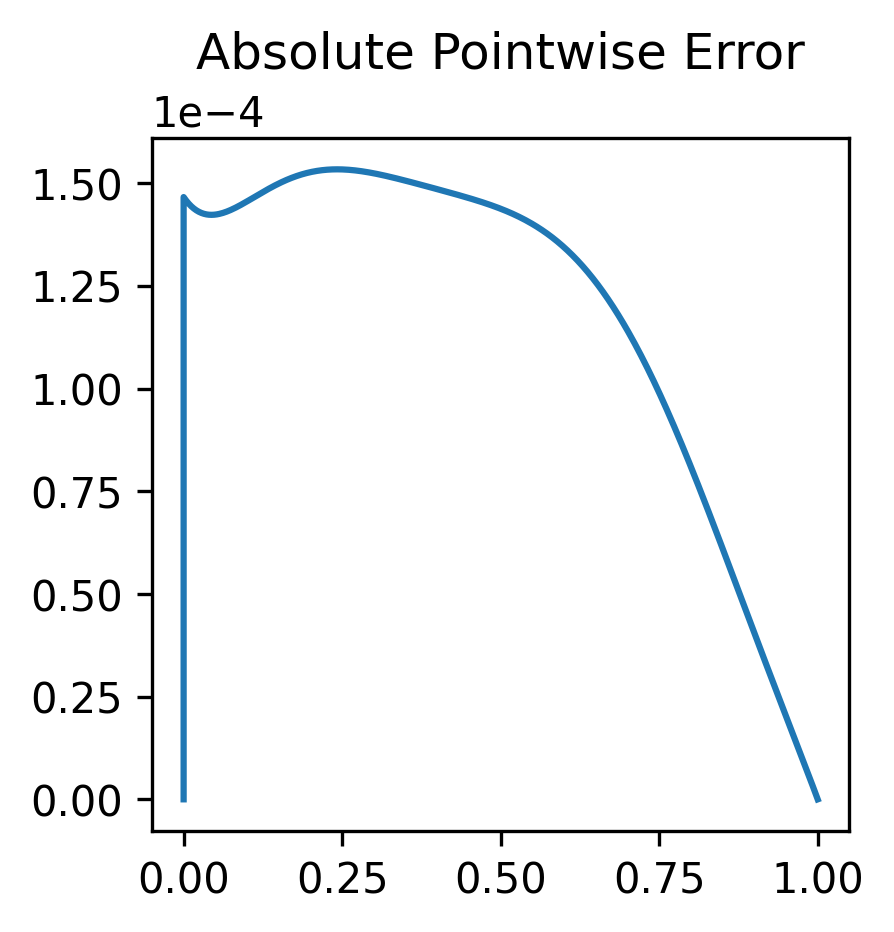}}
\caption{Exact solutions and sl-FVNN predictions of CD eqn.}\label{f:cd_slFVNN}
\end{figure}

\begin{table}[]
\begin{tabular}{|l|ll|ll|}
\hline
                      & FVNNs                                       &                             & {sl-FVNNs}                                    &                             \\ \hline
                      & \multicolumn{1}{l|}{Relative $L^{2}$ error} & Relative $L^{\infty}$ error & \multicolumn{1}{l|}{Relative $L^{2}$ error} & Relative $L^{\infty}$ error \\ \hline 
$\varepsilon=10^{-4}$ & 
\multicolumn{1}{l|}
{4.7539e+00}             & 1.4499e+01                  & \multicolumn{1}{l|}
{1.6323e-04}             & 2.0943e-04                  \\
$\varepsilon=10^{-6}$ & 
\multicolumn{1}{l|}
{1.1403e+02}              & 4.1133e+02                  & \multicolumn{1}{l|}
{1.5543e-04}             & 1.7654e-04                  \\
$\varepsilon=10^{-8}$ & 
\multicolumn{1}{l|}
{7.5334e+02}              & 3.1624e+03                  & \multicolumn{1}{l|}
{1.5718e-04}             & 1.5673e-04                  \\
$\varepsilon=10^{-10}$ & 
\multicolumn{1}{l|}
{1.1918e+04}              & 3.8459e+04                  & \multicolumn{1}{l|}
{1.5422e-04}             & 1.8800e-04                  \\
$\varepsilon=10^{-12}$ & 
\multicolumn{1}{l|}
{9.2081e+04}              & 3.1811e+05                  & \multicolumn{1}{l|}
{1.5632e-04}             & 1.5335e-04                  \\
\hline
\end{tabular}
\caption{Error comparison between FVNNs and {sl-FVNNs} for CD eqn.}
\label{table1_cd}
\end{table}

\begin{table}[]
\begin{tabular}{|l|l|l|}
\hline
                      & 
                      PINNs (Relative $L^{2}$ error)
                      & 
                      {sl-PINNs} (Relative $L^{2}$ error)
                      \\ \hline
                      $\ep = 10^{-4}$
                      & 1.097e+00
                      & 5.474e-02 
                      \\ \hline  
\end{tabular}
\caption{
{
Error comparison between PINNs and {sl-PINNs} for  CD eqn. cited from \cite{GHJ}.
}
}
\label{table1_cd_PINN}
\end{table}

\begin{table}[]
\begin{tabular}{|l|ll|ll|}
\hline
                      & FVNNs                                       &                             & {sl-FVNNs}                                    &                             \\ \hline
                      & \multicolumn{1}{l|}{Training time} & Iterations/sec & \multicolumn{1}{l|}{Training time} & Iterations/sec \\ \hline 
$\varepsilon=10^{-4}$ & 
\multicolumn{1}{l|}
{1.61 sec}             & 1814.10                  & \multicolumn{1}{l|}
{0.71 sec}             & 1300.89                  \\
$\varepsilon=10^{-6}$ & 
\multicolumn{1}{l|}
{1.71 sec}              & 4288.24                  & \multicolumn{1}{l|}
{0.90 sec}             & 1576.65                  \\
$\varepsilon=10^{-8}$ & 
\multicolumn{1}{l|}
{6.95 sec}              & 978.71                  & \multicolumn{1}{l|}
{0.70 sec}             & 3944.69                  \\
$\varepsilon=10^{-10}$ & 
\multicolumn{1}{l|}
{26.64 sec}              & 693.95                  & \multicolumn{1}{l|}
{0.82 sec}             & 514.26                  \\
$\varepsilon=10^{-12}$ & 
\multicolumn{1}{l|}
{13.47 sec}              & 4500.07                  & \multicolumn{1}{l|}
{0.85 sec}             & 502.78                  \\
\hline
\end{tabular}
\caption{Training time comparison between FVNNs and {sl-FVNNs} for CD eqn.}
\label{table2_cd}
\end{table}

\newpage

\subsection{Stationary Burgers Equation (Burgers eqn)}\label{S.comp_Burgers}
\noindent\\
For the computations of (\ref{e:Burgers}), 
we introduce an exact solution $u_{exact}$  in the form, 
\begin{equation}\label{exact_Burgers}
    u_{exact}(x) 
        = 
        - 1 
        - (
        e^{- x/\ep} - x \, e^{- 1/\ep} 
        ), 
\end{equation}
where the data $f$ is properly chosen by computing the right-hand side of (\ref{e:Burgers}) with $u = u_{exact}$.

The training loss values of FVNN and sl-FVNN  at $\ep = 10^{-8}$ appear in 
Figure \ref{f:loss_Burgers} and   
FVNN predictions and sl-FVNN predictions at different $\ep$ values are plotted in Figures \ref{f:Burgers_FVNN} and \ref{f:Burgers_slFVNN}. 
The comparison of numerical errors and training times between FVNNs and sl-FVNNs are summarized in Tables \ref{table1_Burgers} and \ref{table2_Burgers} respectively. 

The Burgers equation presents significant challenges due to steep boundary layers near $x = 0$, especially in the small $\ep$ regime. 
As shown in Figure \ref{f:loss_Burgers}, standard FVNNs exhibit stagnation in loss over tens of thousands of iterations and fail to converge meaningfully. 
In contrast, sl-FVNNs rapidly achieve machine-level accuracy in under $50$ iterations.

Figures \ref{f:Burgers_FVNN} and \ref{f:Burgers_slFVNN}  show that FVNNs produce unstable and inaccurate predictions, with pointwise errors reaching 
$10^2 - 10^3$, especially near $x=0$. 
sl-FVNNs, however, deliver highly accurate approximations across all $\ep$ values, with smooth solutions and uniform error below 
$
10^{-6} 
$.

As summarized in Tables \ref{table1_Burgers} and \ref{table2_Burgers}, sl-FVNNs achieve $4 - 6$ orders of magnitude smaller errors and up to $500$ times faster training time. 
These results confirm that sl-FVNNs are exceptionally well-suited for solving Burgers-type singular problems with steep gradients.

From Tables \ref{table1_Burgers} and \ref{table1_Burgers_PINN}, we see that the sl-FVNN and sl-PINN methods yield comparable accuracy at $\ep = 10^{-8}$, both achieving relative $L^2$ errors on the order of $10^{-3}$, showing their effectiveness in resolving sharp layers over the PINNs and FVNNs.

\begin{figure}
\centering
\subfloat[FVNNs]{\includegraphics[width=0.5\linewidth]{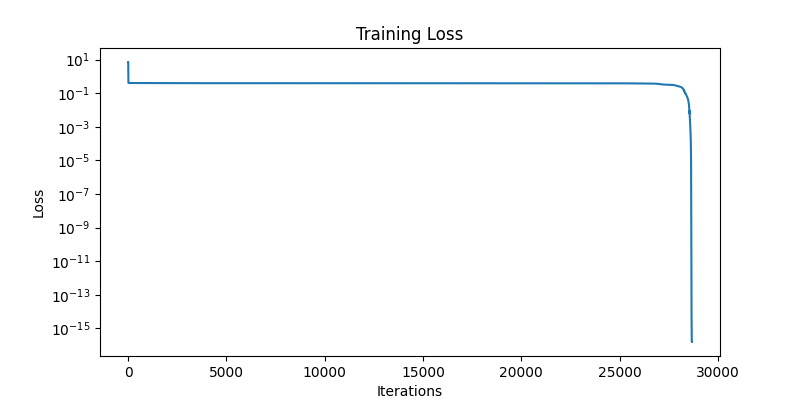}}
\hfill
\subfloat[sl-FVNNs]{\includegraphics[width=0.5\linewidth]{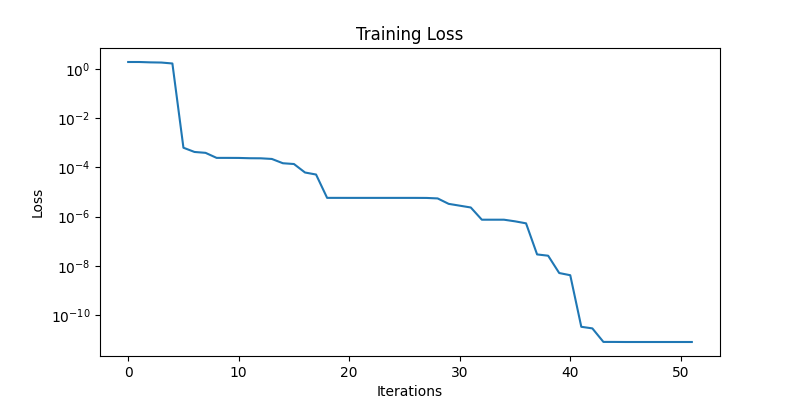}}
\caption{
{
Loss values during the training of Burgers eqn. at $\ep = 10^{-8}$.
}
}\label{f:loss_Burgers}
\end{figure}

\begin{figure}
\centering
\subfloat[$\ep = 10^{-4}$]{\includegraphics[width=0.3\linewidth]{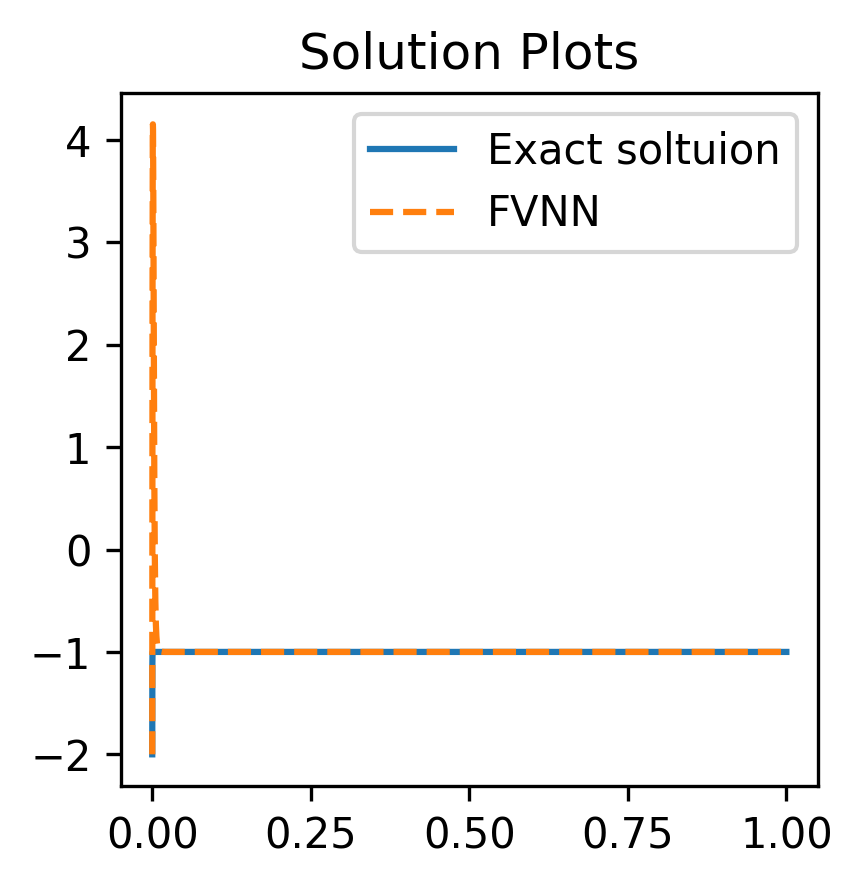}}
\subfloat[$\ep = 10^{-8}$]{\includegraphics[width=0.3\linewidth]{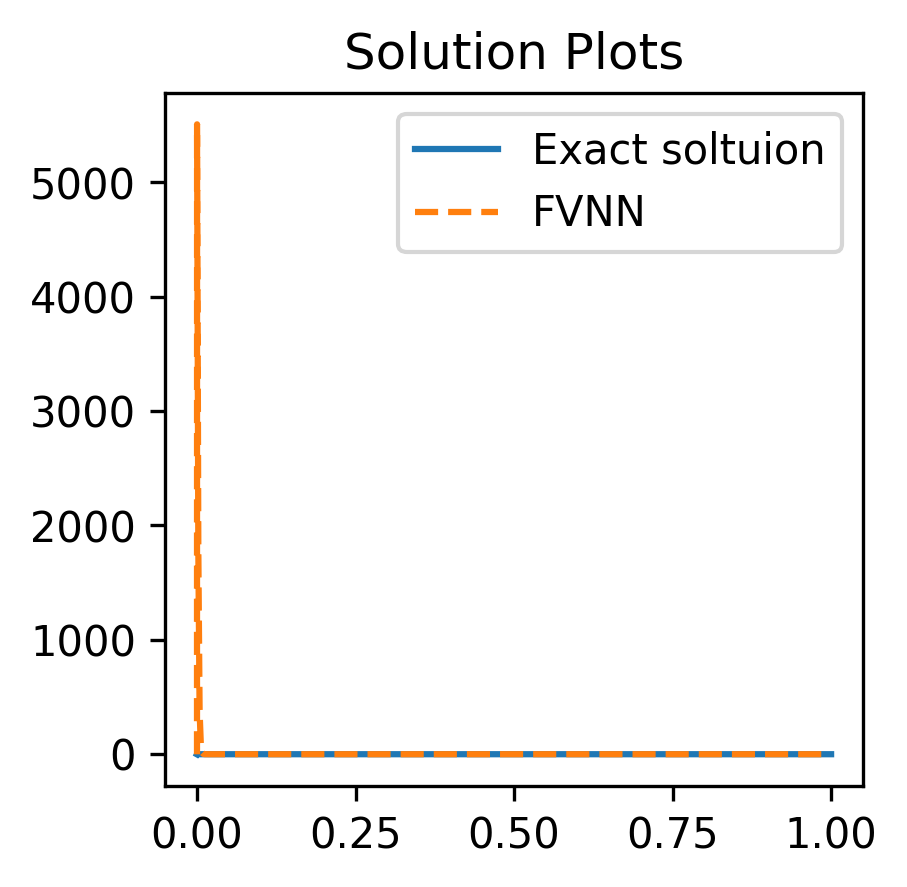}}
\subfloat[$\ep = 10^{-12}$]{\includegraphics[width=0.3\linewidth]{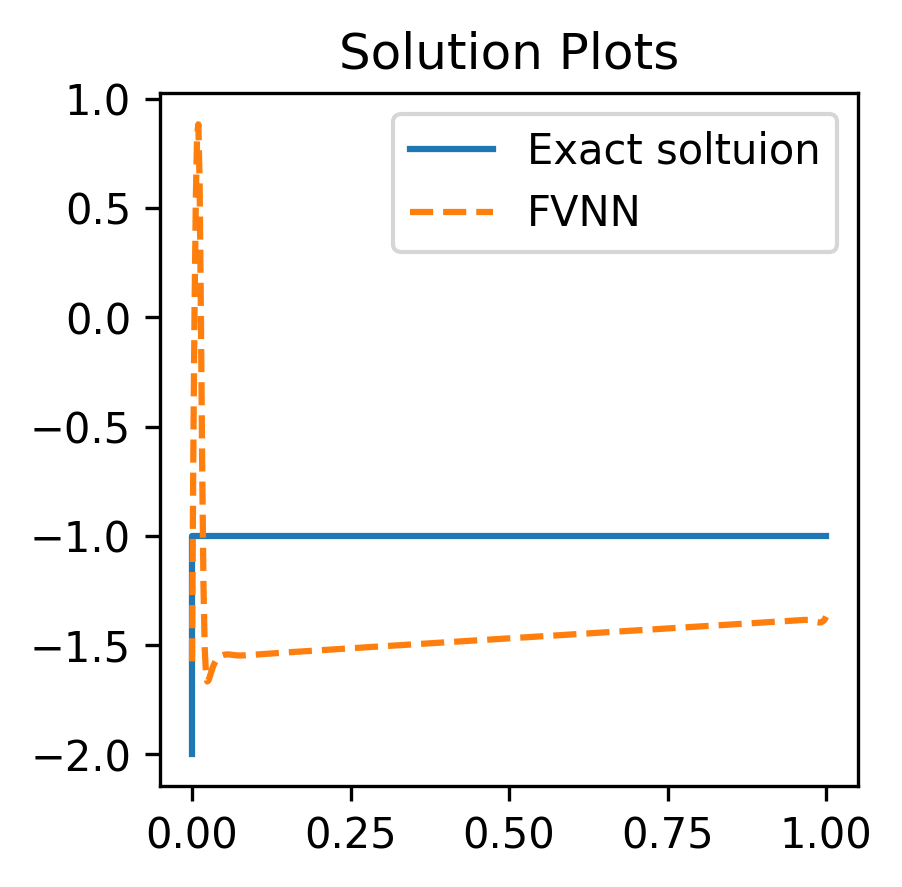}}
\vfill
\subfloat[$\ep = 10^{-4}$]{\includegraphics[width=0.3\linewidth]{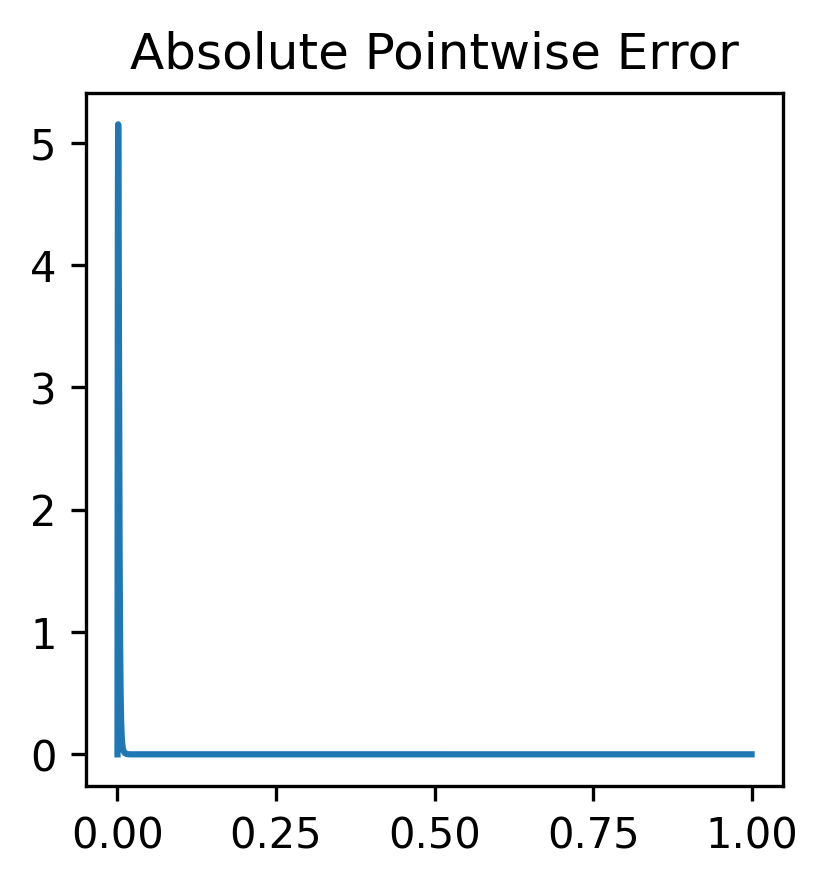}}
\subfloat[$\ep = 10^{-8}$]{\includegraphics[width=0.3\linewidth]{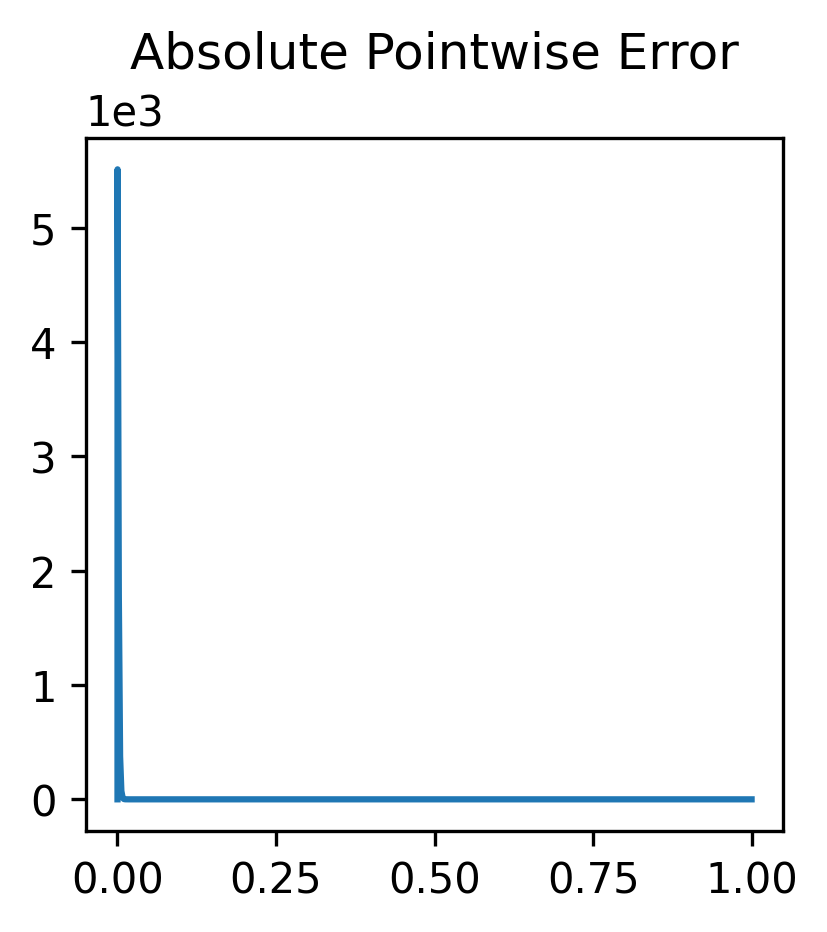}}
\subfloat[$\ep = 10^{-12}$]{\includegraphics[width=0.3\linewidth]{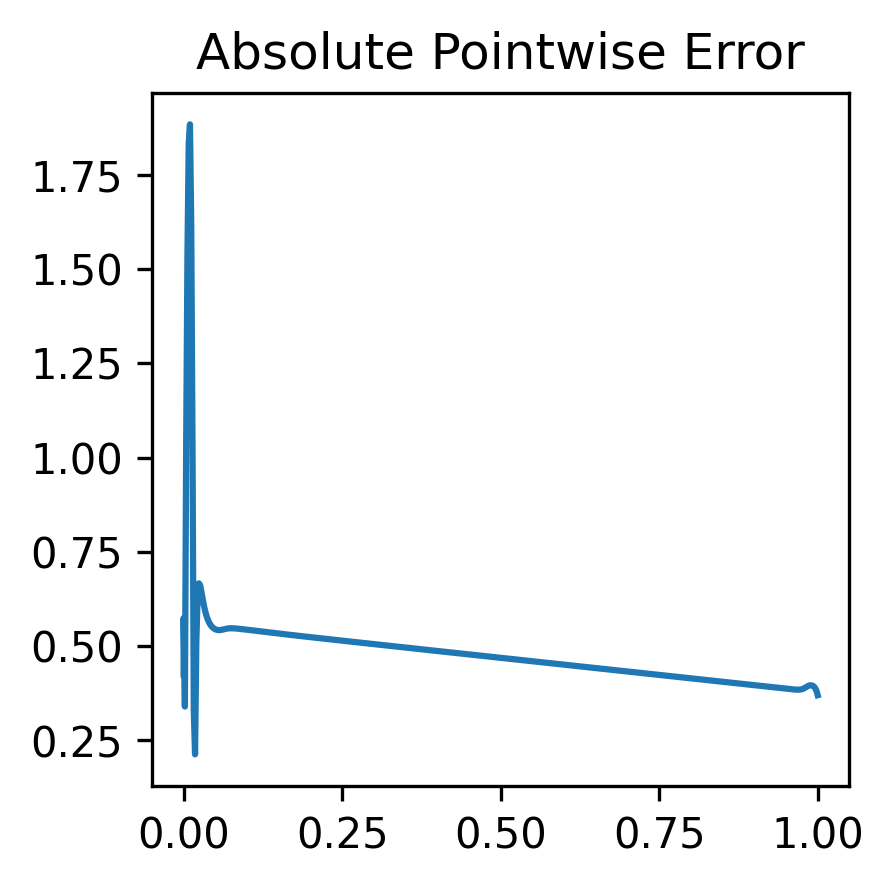}}
\caption{Exact solutions and FVNN predictions of Burgers eqn.}\label{f:Burgers_FVNN}
\end{figure}

\begin{figure}
\centering
\subfloat[$\ep = 10^{-4}$]{\includegraphics[width=0.3\linewidth]{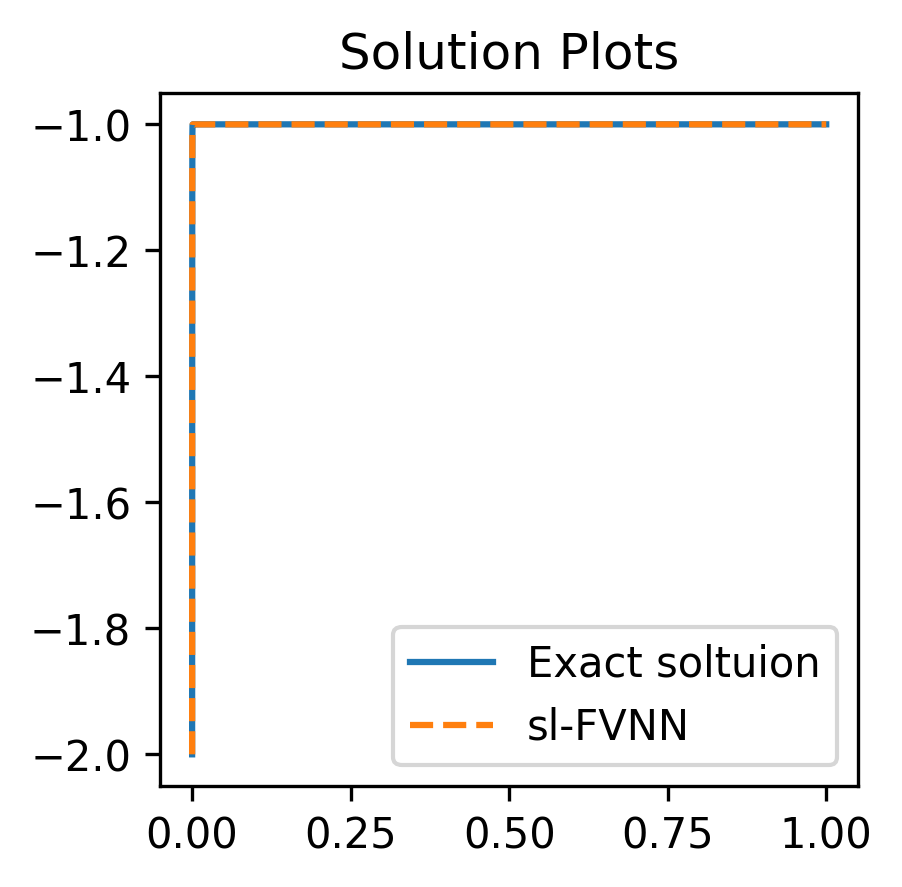}}
\subfloat[$\ep = 10^{-8}$]{\includegraphics[width=0.3\linewidth]{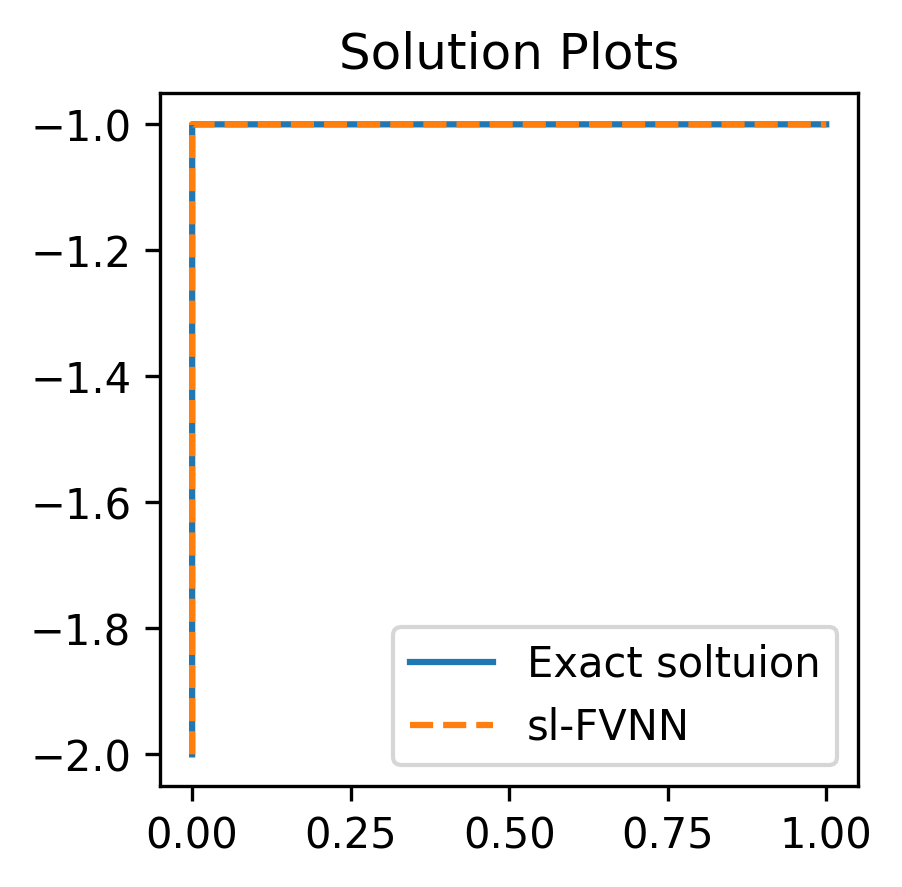}}
\subfloat[$\ep = 10^{-12}$]{\includegraphics[width=0.3\linewidth]{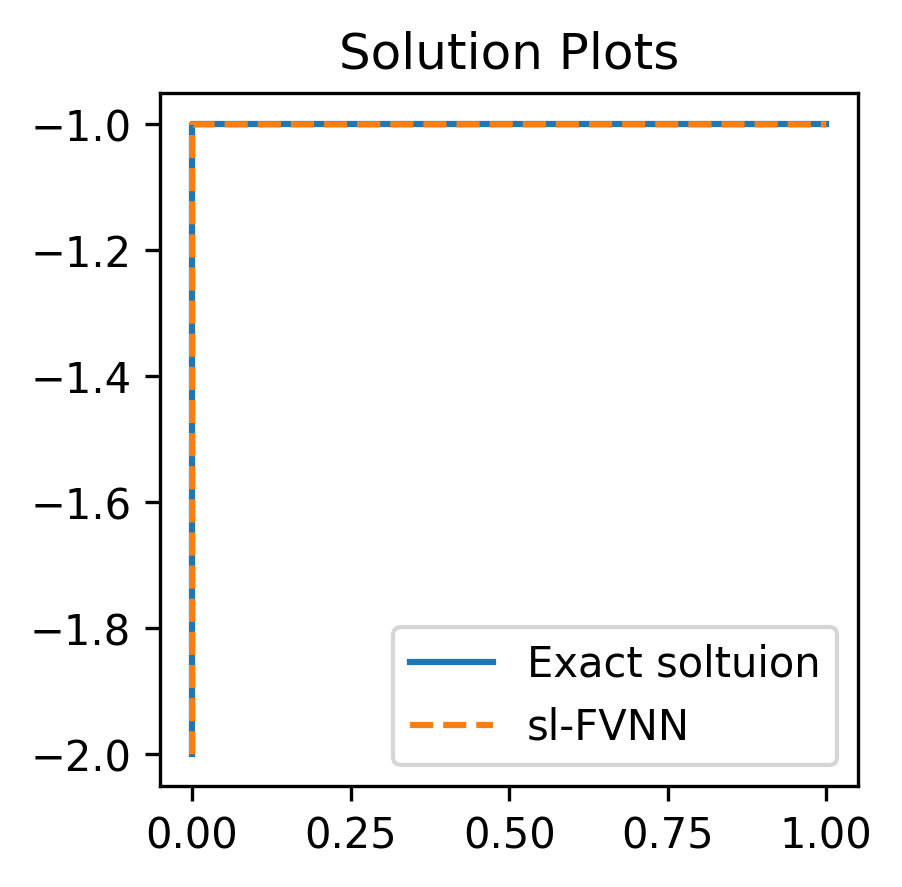}}
\vfill
\subfloat[$\ep = 10^{-4}$]{\includegraphics[width=0.3\linewidth]{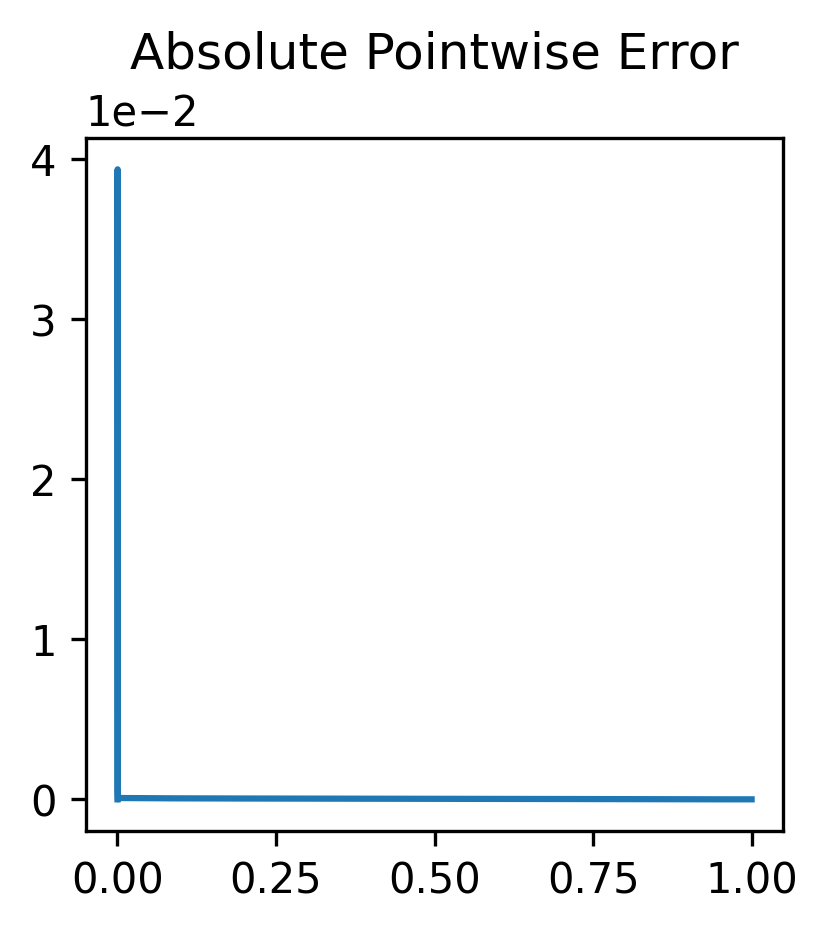}}
\subfloat[$\ep = 10^{-8}$]{\includegraphics[width=0.3\linewidth]{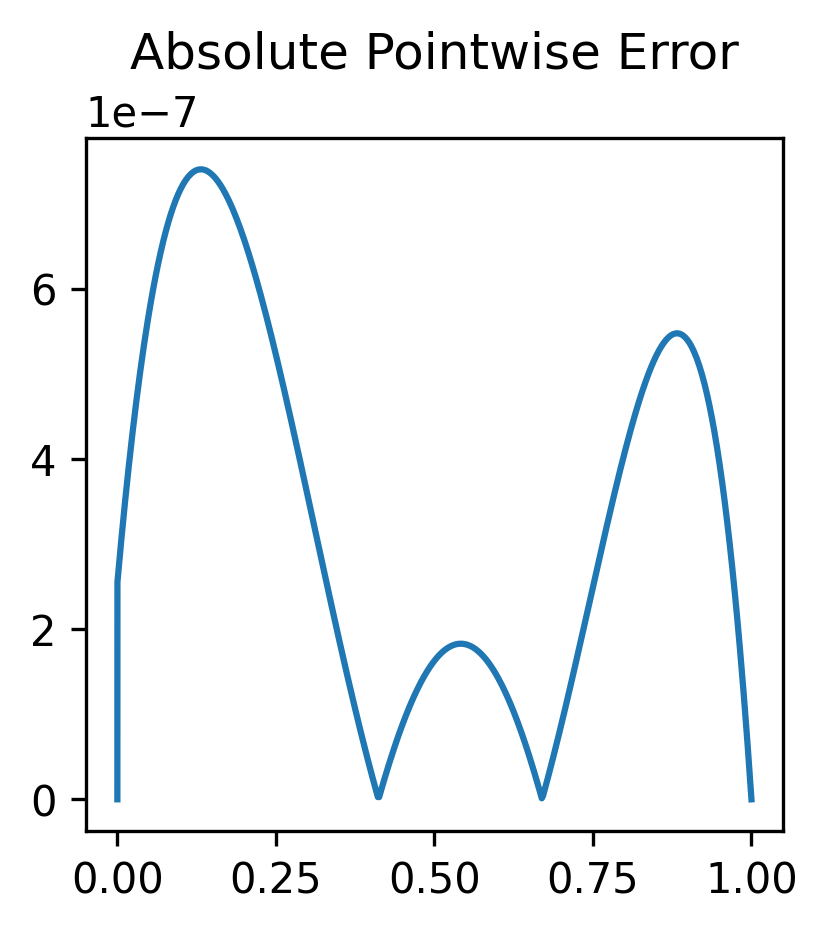}}
\subfloat[$\ep = 10^{-12}$]{\includegraphics[width=0.3\linewidth]{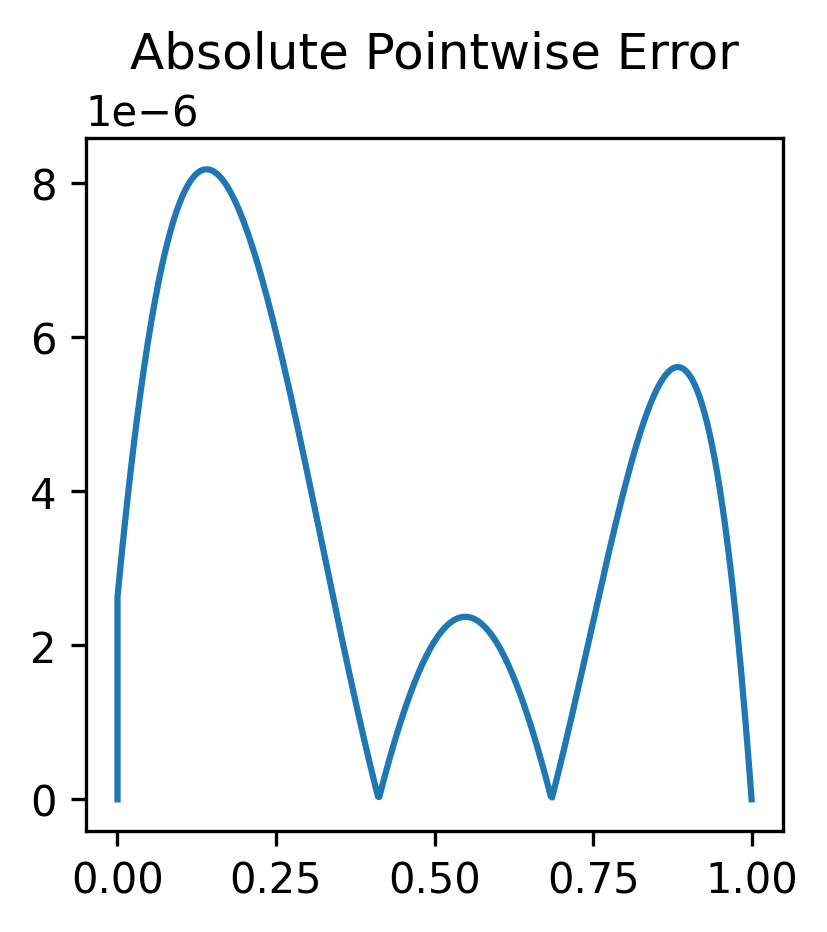}}
\caption{
{
Exact solutions and sl-FVNN predictions of Burgers eqn.
}}\label{f:Burgers_slFVNN}
\end{figure}

\begin{table}[]
\begin{tabular}{|l|ll|ll|}
\hline
                      & FVNNs                                       &                             & {sl-FVNNs}                                    &                             \\ \hline
                      & \multicolumn{1}{l|}{Relative $L^{2}$ error} & Relative $L^{\infty}$ error & \multicolumn{1}{l|}{Relative $L^{2}$ error} & Relative $L^{\infty}$ error \\ \hline 
$\varepsilon=10^{-4}$ & 
\multicolumn{1}{l|}
{6.2117e-01}             & 2.5760e+00                  & \multicolumn{1}{l|}
{1.6204e-03}             & 1.9667e-02                  \\
$\varepsilon=10^{-6}$ & 
\multicolumn{1}{l|}
{1.3843e+02}              & 3.0794e+02                  & \multicolumn{1}{l|}
{4.2945e-06}             & 4.0148e-06                  \\
$\varepsilon=10^{-8}$ & 
\multicolumn{1}{l|}
{1.2409e+03}              & 2.7548e+03                  & \multicolumn{1}{l|}
{3.7953e-07}             & 3.7005e-07                  \\
$\varepsilon=10^{-10}$ & 
\multicolumn{1}{l|}
{7.3086e+01}              & 3.6492e+02                  & \multicolumn{1}{l|}
{3.9992e-06}             & 4.1146e-06                  \\
$\varepsilon=10^{-12}$ & 
\multicolumn{1}{l|}
{4.9958e-01}              & 9.4170e-01                  & \multicolumn{1}{l|}
{4.2479e-06}             & 4.0912e-06                  \\
\hline
\end{tabular}
\caption{Error comparison between FVNNs and {sl-FVNNs} for Burgers eqn.}
\label{table1_Burgers}
\end{table}

\begin{table}[]
\begin{tabular}{|l|l|l|}
\hline
                      & 
                      PINNs (Relative $L^{2}$ error)
                      & 
                      {sl-PINNs} (Relative $L^{2}$ error)
                      \\ \hline
                      $\ep = 10^{-4}$
                      & 1.003e+00
                      & 2.599e-03 
                      \\ \hline  
\end{tabular}
\caption{
{
Error comparison between PINNs and {sl-PINNs} for  Burgers eqn. cited from \cite{GHJ}.
}}
\label{table1_Burgers_PINN}
\end{table}

\begin{table}[]
\begin{tabular}{|l|ll|ll|}
\hline
                      & FVNNs                                       &                             & {sl-FVNNs}                                    &                             \\ \hline
                      & \multicolumn{1}{l|}{Training time} & Iterations/sec & \multicolumn{1}{l|}{Training time} & Iterations/sec \\ \hline 
$\varepsilon=10^{-4}$ & 
\multicolumn{1}{l|}
{1.57 sec}             & 1501.22                  & \multicolumn{1}{l|}
{0.14 sec}             & 380.97                  \\
$\varepsilon=10^{-6}$ & 
\multicolumn{1}{l|}
{26.10 sec}              & 335.32                  & \multicolumn{1}{l|}
{0.05 sec}             & 947.66                  \\
$\varepsilon=10^{-8}$ & 
\multicolumn{1}{l|}
{78.86 sec}              & 363.64                  & \multicolumn{1}{l|}
{0.15 sec}             & 344.33                  \\
$\varepsilon=10^{-10}$ & 
\multicolumn{1}{l|}
{9.53 sec}              & 2822.45                  & \multicolumn{1}{l|}
{0.14 sec}             & 641.04                  \\
$\varepsilon=10^{-12}$ & 
\multicolumn{1}{l|}
{19.34 sec}              & 197.94                  & \multicolumn{1}{l|}
{0.10 sec}             & 911.61                  \\
\hline
\end{tabular}
\caption{Training time comparison between FVNNs and {sl-FVNNs} for Burgers eqn.}
\label{table2_Burgers}
\end{table}

\section{Conclusion}\label{S.Con}
In this article, we developed a singular layer enhanced Finite Volume (Physics Informed) Neural Network (sl-FVNN) framework for solving singularly perturbed boundary value problems. 
By augmenting the standard Finite Volume Neural Network (FVNN) architecture with boundary layer correctors derived from asymptotic analysis, our new sl-FVNN method achieves substantial improvements in both accuracy and efficiency across a wide range of small perturbation  parameters $\ep$. 
The sl-FVNN consistently outperforms traditional FVNNs, particularly in stiff regimes when $\ep$ is small, as demonstrated by orders of magnitude reductions in both relative error and training time. 

Our comparative studies further show that both sl-FVNNs and sl-PINNs (studied in \cite{GHJ}) perform well on benchmark problems such as the reaction-diffusion and Burgers equations, attaining relative errors of the same order. However, for convection-dominated nonlinear problems, sl-FVNNs significantly outperform sl-PINNs, achieving approximately $10^2$ times smaller relative $L^2$ error. These results confirm that embedding singular layer correctors into neural network formulations provides a powerful and efficient methodology for solving multiscale PDEs without requiring excessively deep or over parameterized architectures. 

Future work will extend this framework to higher-dimensional and time-dependent problems in complex domains. 

{
We also note that, for the one-dimensional benchmark problems considered in this article, once the leading boundary layer has been analytically identified, classical numerical methods, e.g., finite difference or finite volume methods 
can also be effective for the remaining smooth problem. 
The purpose of the present work is not to replace such classical schemes in simple one-dimensional settings, but 
to show how asymptotic layer correctors can be embedded into a Finite Volume Neural Network (FVNN) framework that incorporates singular-layer information and remains stable in small-$\ep$ regimes.
}

\section*{Acknowledgment}
{
The authors would like to thank the reviewers for their valuable comments and suggestions, which helped improve the quality of the manuscript. 
The first author was partially supported by the Fulbright U.S. Scholar Program, which is sponsored by the U.S. Department of State and the Korean-American Educational Commission. 
The first author also gratefully acknowledges the hospitality of Seoul National University during his visit.
}
 
\section*{Availability of data}
The data that support the findings of this study are available from the corresponding author, Dr. Junghwa Kim, upon reasonable request.

\end{document}